\documentclass[preprint,12pt]{elsarticle}

\usepackage{graphicx}
\usepackage{epsfig}
\usepackage{tikz}

\usepackage{amssymb}
\usepackage{amsthm}

\usepackage{lineno}
\usepackage[english]{inputenc}
\usepackage{fancyhdr}
\usepackage{indentfirst}
\usepackage{color}
\usepackage{newlfont}
\usepackage{textcomp}
\usepackage{amsmath}
\usepackage{latexsym}
\usepackage{mathrsfs}
\usepackage{bbold}
\usepackage{listings}
\usepackage{caption}
\usepackage{subcaption}

\usepackage[subnum]{cases}
\usepackage{calc}

\usepackage{floatrow}
\usepackage{enumitem}

\usepackage{algpseudocode,algorithm}

\usepackage{lipsum}
\makeatletter
\def\ps@pprintTitle{%
 \let\@oddhead\@empty
 \let\@evenhead\@empty
 \def\@oddfoot{}%
 \let\@evenfoot\@oddfoot}
\makeatother

\newcommand{\footnoteremember}[2]{%
  \footnote{#2}%
  \global\expandafter\edef\csname #1@fnnumber\endcsname{\thefootnote}%
}

\newcommand{\footnoterecall}[1]{%
  \footnotemark[\csname #1@fnnumber\endcsname]%
}

\usetikzlibrary{calc,trees,positioning,arrows,chains,shapes.geometric,%
    decorations.pathreplacing,decorations.pathmorphing,shapes,%
    matrix,shapes.symbols}

\tikzset{
>=stealth',
  punktchain/.style={
    rectangle, 
    rounded corners, 
    draw=black, very thick,
    text width=10em, 
    minimum height=3em, 
    text centered, 
    on chain},
  line/.style={draw, thick, <-},
  element/.style={
    tape,
    top color=white,
    bottom color=blue!50!black!60!,
    minimum width=10em,
    draw=blue!40!black!90, very thick,
    text width=10em, 
    minimum height=3.5em, 
    text centered, 
    on chain},
  every join/.style={->, thick,shorten >=1pt},
  decoration={brace},
  tuborg/.style={decorate},
  tubnode/.style={midway, right=2pt},
}

\newtheorem{theorem}{Theorem}

\theoremstyle{remark}

\begin{document}





\begin{frontmatter}
\title{A highly accurate boundary integral equation method for surfactant-laden drops in 3D}

\author[kth]{Chiara Sorgentone}
\ead{sorgento@kth.se}

\author[kth]{Anna-Karin Tornberg}
\ead{akto@kth.se}

\address[kth]{KTH Mathematics, Linn\'e Flow Centre/Swedish e-Science Research Centre, 10044 Stockholm, Sweden}
	

 \begin{abstract} 
The presence of surfactants alters the dynamics of viscous drops immersed in an ambient viscous fluid. This is specifically true at small scales, such as in applications of droplet based microfluidics, where the interface dynamics become of increased importance. At such small scales, viscous forces dominate and inertial effects are often negligible.  Considering Stokes flow, a numerical method based on a boundary integral formulation is presented for simulating 3D drops covered by an insoluble surfactant. The method is able to simulate drops with different viscosities and close interactions, 
automatically controlling the time step size and maintaining high accuracy also when substantial drop deformation appears. To achieve this, the drop surfaces as well as the surfactant concentration on each surface are represented by spherical harmonics expansions. A novel reparameterization method is introduced to ensure a high-quality representation of the drops also under deformation, specialized quadrature methods for singular and nearly singular integrals that appear in the formulation are evoked and the adaptive time stepping scheme for the coupled drop and surfactant evolution is designed with a preconditioned implicit treatment of the surfactant diffusion.

 \end{abstract}

\end{frontmatter}



\section*{Introduction}
\label{Introduction}

Micro-fluidics, or fluid mechanics at small scales, is an expanding
research area.  Droplet-based microfluidics is a subcategory of
microfluidics where tiny (e.g. picoliter sized) droplets are used as
the equivalent of ``test tubes'', offering the possibility of high
troughput experiments that enhance the speed of chemical and
biological assays \cite{microfluidics,microfluidics_new}.
Owing to the small scales, viscous forces dominate and intertial
effects are often negligible. Due to the large surface to volume
ratio, the importance of interface dynamics becomes pronounced. As
water droplets are dispersed in oil, they are stabilized by a powerful
surfactant preventing them from coalescing. These surface active
agents are compounds that lower the surface tension between liquids;
they are widely used in engineering applications, in pharmaceuticals,
foods and petroleum industries \cite{surf1,surf2}. Surfactants
can also be found in e.g. detergents, emulsifiers, paints, adhesives,
inks and alveoli.

The surfactant concentration can be modeled with a time dependent differential equation defined on the moving and deforming interface, coupled with the Stokes equation for the evolution of the droplets through a non-linear equation of state relating the surface tension to the surfactant concentration. 
There are two competing processes governing the effect of surfactants on drop deformation \cite{feigl}: surfactant dilution and surfactant convection. The first one is due to deformation of the drop from a spherical shape, which causes an increase of total surface area and a consequent decrease of surfactant concentration, hence smaller deformation. The fluid flow will however typically cause surfactant convection that tends to increase the concentration at the tips of the drop, resulting in an increasing drop deformation.

There are several studies of drop behavior without surfactants (see \cite{drops_review,drops_bubble_book,summary_microfluidics} for a summary), but the behavior of the surfactant-covered droplets has been much less studied, especially in 3D.
In 2D there is a number of papers for insoluble surfactants, see
e.g. \cite{Stone2D,Xu2D,AK_2D_surfactant}. For soluble
surfactants where the surfactant is assumed to have a concentration
also in one of the fluid phases, there are fewer studies, see
e.g. \cite{AK-2D-soluble,Xu-2D-soluble}. 
For the 3D problem, one of the earlier studies is due to Li and
Pozrikidis \cite{lipoz} in 1997, followed by another paper of Yon and
Pozrikidis \cite{yonpoz}, where the drop deformation in shear flow is
studied only for the special case of viscosity ratio equal to one
(same viscosity inside and outside the drop). Arbitrary viscosity
ratio has been considered later on:  Bazhlekov et al. presented a
numerical investigation of the effect of insoluble surfactants on drop
deformation and break up \cite{18} and in \cite{feigl} Feigl et al. presented a
study of the same problem combining numerical simulations with
experiments on droplets. In all of the mentioned cases a
boundary integral formulation has been used. More recently, Teigen et
al. proposed a diffuse-interface method \cite{soluble_3D}, and Muraduglo and Tryggvason \cite{Muradoglu} a
front-tracking method for drops with soluble surfactant. 

These micro-fluidic problems are very challenging to simulate accurately and efficiently, especially in three dimensions. 
The geometry is time-dependent due to the motion and deformation of the drops, surface tension forces induce jumps in pressure and velocity gradients across the drop surfaces, and the presence of surfactants require the solution of a time dependent PDE on each deforming drop surface. 
Using grid-based methods to solve the Stokes (or Navier-Stokes)
equations, the computational grid/mesh is typically not required to
conform to the the drop interfaces that could be kept track of by
front-tracking \cite{Muradoglu}, a level set method \cite{Xu2D} or a phase field function with a diffuse interface layer as in \cite{soluble_3D}.  The singular surface tension forces are commonly regularized onto the grid, yielding only first order accuracy close to the interface as measured in the resolution of the volume grid. Recent advances has been made to develop so-called cut cell methods to increase this accuracy to second order, see e.g. \cite{cut_zahedi} and the references therein. 

For the case of Stokes equations a boundary integral equation can be
formulated. This is a sharp interface formulation
for which the above mentioned issues are avoided - there is no
underlying volume grid to couple to, jumps in solutions are naturally
taken care of, and viscosity ratios between fluids enter only in
coefficients of the equations. The formulation contains integrals only
over the fluid-fluid interfaces and possibly external boundaries,
reducing the dimension of the problem, and hence the number of
unknowns in the discretized problem. 

There are many discretization choices to be made for a boundary
integral formulation, including how the drop surfaces should be
represented and how the integrals in the equation should be evaluated.
The integrals contain Green's functions for Stokes equations, which
are singular but integrable when integrating over a drop surface for
an evaluation point at that same surface, and nearly singular for an
evaluation point on a drop surface nearby. Evaluating such integrals
require special techniques.

One common way to discretize a surface is to use a triangulated
surface \cite{yonpoz,zinchenko}. We have chosen to use a
higher order representation in terms of a parameterization of the
surface, expressed as a spherical harmonics expansion with vector
valued coefficients for the three drop coordinates. This has earlier
been used e.g for the simulation of vesicles and blood cells in Stokes
flow \cite{1,6,21}.  An expansion of order $p$
corresponds to $O(p^2)$ discrete points on the surface, and a discrete
spherical harmonics transform can be used to compute the coefficients
from the coordinates of these points.  To evaluate the singular
integrals mentioned above, we use a method that yields exponential
convergence in $p$ \cite{9,6}.  For the case of close interactions and
hence nearly singular integrals, we adapt and optimize a procedure
first presented by Ying et al. \cite{13} to our problem. The errors
can be kept at a very low level also here, if parameters are chosen
appropriately.

When a drop moves and deforms, the point distribution will in general
be distorted, and the surface needs to be reparametrized. The optimal
spherical harmonics representation is the one that describes the
surface using a set of expansion coefficients with the fastest decay
with $p$, corresponding to a specific set of discrete points on the
surface.  The reparameterization procedure should yield a distribution
close to optimal, and conserve the volume and the area of the drop.
In addition, it needs to yield a new parameterization for the
surfactant concentration, that is also represented using a spherical
harmonics expansion. The reparameterization procedure introduced by Veerapaneni et al. in \cite{6} sets
a nice framework which we have used to design a novel
reparameterization technique that meets our requirements regarding
volume/area conservation and reparametrization of the surfactant
representation.

An adaptive time-stepping scheme is used to couple together the drop
and the surfactant evolution. It uses an implicit treatment for the
diffusive term in the surfactant equation, which reduces the CFL
condition to first order.  To make this implicit solver efficient also
for large diffusion coefficients, a new preconditioner is introduced.

In total, the Boundary Integral Method (BIM) for 3D drops covered by an
insoluble surfactant that we present here is able to handle drops with
different viscosities and close interactions maintaining high accuracy
even when strong deformations appears, automatically controlling the
time stepping size.  

The paper is organized as follows: in Section
\ref{math_formulation} we introduce the model, the boundary integral
formulation and the Galerkin method used to solve the discrete
system. In Section \ref{num_int} we discuss the numerical integration
with special attention to singular and nearly singular integrals.  A
study and comparison of different timestepping schemes coupling
together the drop and the surfactant evolution is presented in Section
\ref{time}, where we also introduce the new preconditioner for the
surfactant equation. The new reparametrization procedure is explained
in Section \ref{reparametrization} and numerical experiments for
both clean and surfactant-covered drops are presented in Section
\ref{experiments}, where we validate our method against numerical and
experimental results and where we also show the effect of surfactant
in very near interaction of two drops.  Conclusions and future
perspectives are given in the last section.

\section{Mathematical formulation}
\label{math_formulation}

We consider $N$ drops immersed in a liquid with
different viscosity. The Stokes equations read 
\begin{equation}
\label{eq:diff_eq1}
\begin{cases}
-\mu_i \Delta \textbf{u}(\textbf{x})+\nabla P(\textbf{x}) =0 \\ 
\nabla \cdot \textbf{u}(\textbf{x})=0 \\
\textbf{u}(\textbf{x}) \rightarrow \textbf{u}_\infty(\textbf{x}) \text{ as } |\textbf{x}|\rightarrow \infty
\end{cases}
\end{equation}
for every $\textbf{x}$ inside the $i$-th drop ($i=1,\dots,N$) or in the exterior region ($i=0$), where $\textbf{u}$ is the fluid velocity, $P$ is the pressure and $\mu_i$ is the viscosity.

Let $\gamma^*=\bigcup_i \gamma_i$ denote the union of all drop
surfaces/interfaces.
There is a jump condition at $\gamma^*$ that is given by
\begin{equation}
\label{eq:diff_eq2}
[[T\textbf{n}]]=\textbf{f}(\textbf{x})
\end{equation}
where $T=-PI+\mu(\nabla \textbf{u}+\nabla \textbf{u}^T)$ is the Cauchy
stress tensor, $[[\cdot]]$ denotes the jump across the interface and
$\textbf{f}=2\sigma(\textbf{x})K(\textbf{x})\textbf{n}(\textbf{x})-\nabla_\gamma(\sigma)$
is the interfacial force, with $\textbf{n}$ being the outward pointing unit normal.
The evolution of the interfaces is given by
\begin{equation}
\label{eq:diff_eq3}
\frac{d\textbf{x}}{dt}=\textbf{u}(\textbf{x}), \text{ for all } \textbf{x} \in \gamma^*.
\end{equation}

\begin{table}[h!]
\begin{tabular}{r l}
  \hline
  Symbol & Definition \\
    \hline
  $\Delta_\gamma$ & Surface Laplacian \\
  $\nabla_\gamma$ & Surface gradient \\ 
  $\gamma^*$ & Union of all drop surfaces $\bigcup_i \gamma_i$\\
  $\mu_0$ & Viscosity of the ambient fluid \\ 
  $\mu_i$ & Viscosity of the fluid inside the $i$-th drop\\
   $\lambda_i$ & The viscosity contrast $\frac{\mu_i}{\mu}$\\
  $\textbf{u}$ & Interfacial velocity\\
  $\textbf{u}_\infty$ & Far field velocity \\
   $\textbf{u}_\gamma$ & Tangential component of the interfacial velocity \\
    $\textbf{u}_n$ & Normal component of the interfacial velocity \\
  $P$ & Pressure\\
  $\sigma$ & Dimensionless interfacial tension \\
  $\Gamma$ & Dimensionless surfactant concentration \\
  $Pe$ & P\'eclet number \\
  \hline
\end{tabular}
  \caption{List of symbols}
\end{table}

\begin{table}[h!]
\begin{tabular}{r l | r l}
  \hline
  Symbol & Definition & Symbol & Definition \\
    \hline
    E & $\textbf{x}_\phi \cdot \textbf{x}_\phi $&L & $\textbf{x}_{\phi \phi} \cdot \textbf{n}$\\
F & $\textbf{x}_\phi \cdot \textbf{x}_\theta $&M &$\textbf{x}_{\phi \theta} \cdot \textbf{n}$\\
G & $\textbf{x}_\theta \cdot \textbf{x}_\theta$ &N &$\textbf{x}_{\theta \theta} \cdot \textbf{n}$\\
W & $\sqrt{EG-F^2}$ & K &$\frac{LN-M^2}{W^2}$, mean curvature\\
$\textbf{n}$ & $\frac{\textbf{x}_\phi \times \textbf{x}_\theta}{W}$, unit normal & H & $\frac{1}{2}\frac{EN-2FM+GL}{W^2}$, Gaussian curvature\\
  \hline
\end{tabular}
  \caption{List of geometrical quantities}
\end{table}

Letting $a$, $\sigma_0$ and $\mu_0$ be, respectively, the radius of
an undeformed spherical droplet, the interfacial tension of the
clean interface and the viscosity of the ambient fluid, we can
non-dimensionalize by taking the characteristic length to be $a$, the
characteristic velocity to be $U=\frac{\sigma_0}{\mu_0}$ and the
characteristic time $T=\frac{a\mu_0}{\sigma_0}$. Using the
corresponding non-dimensionalized variables, equations
(\ref{eq:diff_eq1})-(\ref{eq:diff_eq2}) can be reformulated using a
standard boundary integral approach \cite{10}. 
For all $\textbf{x}_0 \in \gamma_i \ (i=1,\dots,N)$,
 \begin{equation}
 \label{eq:main_eq}
 \begin{aligned}
 (\lambda_i+1)\textbf{u}(\textbf{x}_0)=2\textbf{u}_\infty(\textbf{x}_0)&-\sum_{j=1}^N  \bigg( \frac{1}{4\pi}\int_{\gamma_j} \textbf{f}(\textbf{x})\cdot G(\textbf{x}_0,\textbf{x})ds(\textbf{x})\bigg) \\
 &+\sum_{j=1}^N \bigg(\frac{\lambda_i-1}{4\pi}\int_{\gamma_j} \textbf{u}(\textbf{x})\cdot T(\textbf{x}_0,\textbf{x})\cdot \textbf{n}(\textbf{x}) ds(\textbf{x})\bigg),
 \end{aligned}
\end{equation}
where $\lambda_i$ denotes the viscosity contrast of the $i$-th drop.
The tensors $G$ and $T$ are the Stokeslet and the Stresslet, 
$$G(\textbf{x}_0,\textbf{x})=I/r+\hat{\textbf{x}}\hat{\textbf{x}}/r^3, \ \ \ T(\textbf{x}_0,\textbf{x})=-6\hat{\textbf{x}}\hat{\textbf{x}}\hat{\textbf{x}}/r^5, $$
where $\hat{\textbf{x}}=\textbf{x}_0-\textbf{x}$ and
$r=|\hat{\textbf{x}}|$. 
The limit to the boundary has been taken in (\ref{eq:main_eq}) and the integrals are to be understood as the principal value integrals. Using
Eq. (\ref{eq:main_eq}), we solve for the velocity of the drop
surfaces, whereafter the drop evolution can be computed.
The velocity can be evaluated anywhere in the domain whenever it is
wanted (see \cite{10} for details) but is not needed to determine the
drop evolution. The velocity field is continuous across the interfaces
$\gamma_i$, but its gradient and the pressure are not.

The evolution of the drops is affected by the presence of insoluble
surfactant through a modified surface tension force. The equation
governing the evolution of the surfactant concentration $\Gamma$ is a convective-diffusion equation
which can be derived stating the conservation of surfactant mass
\cite{17}; it is given in dimensionless terms by \cite{17,18},
\begin{equation}
\label{eq:surfactant0}
\frac{\partial \Gamma}{\partial t}+\nabla_\gamma \cdot (\Gamma \textbf{u}_\gamma) -\frac{1}{Pe}\nabla_\gamma^2\Gamma+2K(\textbf{x})\Gamma (\textbf{u}\cdot \textbf{n})=0,
\end{equation}
where $Pe=\frac{a \sigma_0}{\mu_0 D_S}$ is the P\'eclet number, with
$D_S$ being the surface diffusivity. The interfacial tension is
related to the surfactant concentration by the equation of
state. Different equations of state can be used \cite{marangoni}, we
will use the Szyszkowski equation (also called Langmuir equation of
state) which in dimensionless form is given by
\begin{equation}
\label{eq:langmuir}
\sigma(\Gamma)=1+E\ln(1-x_s\Gamma)
\end{equation}
where $E$ is the elasticity number and $x_s$ is the surface coverage,
$0 \leq x_s \leq 1$. Equations
(\ref{eq:surfactant0})-(\ref{eq:langmuir}) were non-dimensionalized by
the characteristic velocity $U$, characteristic time $T$ and the
surfactant concentration $\Gamma_{eq}$ for the same amount of surfactant
uniformly distributed over a spherical drop with the same volume \cite{18}.

\subsection{Surface Representation}
Assuming the surface $\gamma_i$ of each drop to be smooth and of
spherical topology, we will use a spherical harmonics expansion to
represent it. Let $\textbf{x}:U \rightarrow \gamma_i$ be a
parametrization of $\gamma_i$, where the domain U is a rectangle
$\{(\theta,\phi): \theta \in [0,\pi], \phi \in [0,2\pi)\}$. A scalar
spherical harmonic function of degree $n$ and order $m$ is given by
\begin{equation}
Y_n^m(\theta,\phi)=\sqrt{\frac{2n+1}{(4\pi)}\frac{(n-m)!}{(n+m)!}}P_n^{m}(\cos(\theta))e^{im\phi}, \ \ \ -n\leq m \leq n,
\end{equation}
where the associated Legendre function is defined by
\begin{equation}
P_n^{m}(x)=(-1)^{m}(1-x^2)^{m/2}\frac{d^m}{dx^m}P_n(x), \ \ \ m \geq 0,
\end{equation}
and
\begin{equation}
\label{eq:pneg}
P_n^{-m}(x)=(-1)^{m}\frac{(n-m)!}{(n+m)!}P_n^m(x), \ \ \ m \geq 0,
\end{equation}
where $P_n$ is the Legendre polynomial of degree $n$,
\begin{equation}
P_n(x)=(2^nn!)^{-1}\frac{d^n}{dx^n}[(x^2-1)^n].
\end{equation}

Let $f:\mathbb{S}^2 \rightarrow \mathbb{R}$ be a square-integrable function. It can be expanded in terms of spherical harmonics as
\begin{equation}
f(\theta,\phi)=\sum_{n=0}^\infty \sum_{m=-n}^n f_n^m Y_n^m(\theta,\phi)
\end{equation}
where the coefficients $f_n^m$ are the moments of the expansion,
\begin{equation}
\label{eq:back_transf}
f_n^m=(f,Y_n^m)=\int_0^\pi \int_0^{2\pi}f(\theta,\phi)\overline{Y_n^m(\theta,\phi)} \sin \theta d\theta d\phi.
\end{equation}
We truncate the spherical harmonic expansion at some degree $p$, 
\begin{equation}
\label{eq:spharm_ex}
f_p(\theta,\phi)= \sum_{n=0}^p \sum_{m=-n}^n f_n^m Y_n^m(\theta,\phi).
\end{equation}
A drop surface will be represented by three such truncated scalar
expansions, 
\begin{equation}
\label{eq:surface_ev}
\textbf{x}(\theta,\phi)=[x(\theta,\phi),y(\theta,\phi),z(\theta,\phi)]=\sum_{n=0}^p \sum_{m=-n}^n [x_n^m,y_n^m,z_n^m]Y_n^m(\theta,\phi)
\end{equation}
where we have omitted the subscript $p$ and also the index of the
drop for simplicity. We will denote by $\textbf{x}_n^m=[x_n^m,y_n^m,z_n^m]$ the spherical harmonics coefficients vector.
The surfactant concentration corresponding to this drop surface is
similarly represented as 
\begin{equation}
\Gamma(\theta,\phi)=\sum_{n=0}^p \sum_{m=-n}^n \Gamma_n^m Y_n^m(\theta,\phi).
\end{equation}
We will compute derivatives of functions using analytical formulas in the spectral domain, see \ref{appendix1}.

\subsection{Galerkin Formulation}
\label{galerkin}

To solve the integral equation (\ref{eq:main_eq}), the evolution equation (\ref{eq:diff_eq3}) and
the surfactant equation (\ref{eq:surfactant0}), we use a Galerkin formulation as
suggested in \cite{1}, where a similar approach is used to simulate the flow of vesicles. Let $(\cdot,\cdot)$ denote the inner product
$(f,Y_{nm})=\int_\mathbb{S^2} f \overline{Y_{nm}}ds$, then the Galerkin
method seeks the solution to the original system by
\begin{equation}
\label{eq:BIM_galerkin}
(\lambda_i+1)(\textbf{u},Y_{nm})= (\textbf{u}^\infty,Y_{nm})+\sum_{j=1}^N [ (\textit{S}_j[\textbf{f}],Y_{nm})+((\lambda_i-1)\textit{D}_j[\textbf{u}],Y_{nm})],
\end{equation}
\begin{equation}
(\frac{\partial \textbf{x}}{\partial t},Y_{nm}) = (\textbf{u},Y_{nm}),
\end{equation}
\begin{equation}
\label{eq:galerkin_surfactant}
(\frac{\partial \Gamma}{\partial t},Y_{nm}) = -(\nabla_\gamma \cdot \Gamma \textbf{u}_\gamma, Y_{nm}) +\frac{1}{Pe}(\nabla_\gamma^2\Gamma,Y_{nm})-(2K\Gamma \textbf{u}\cdot \textbf{n},Y_{nm}).
\end{equation}
This approach exploits the orthogonality and symmetry properties of the spherical harmonics as discussed in \cite{SPHARM} and it results to be cheaper compared to a standard approach, as observed already in \cite{1}.
Indeed, if the original system with physical variables has size $2(p+1)^2$ (the grid points), the Galerkin formulation almost halves the number of unknowns (the spherical coefficients). The spectral coefficients are complex (doubling the size of system), but thanks to the skew-symmetric property of real-valued functions,
\begin{equation*}
f_n^m=(-1)^m\overline{f}_n^{-m}
\end{equation*}
and the fact that all coefficients with $|m|>n$ are zero, the final system has $(p+1)(p+2)$ unknowns.

\subsection{Upsampling and De-aliasing} 
\label{upsampling_dealiasing}

Spherical harmonics can be seen as the surface counterpart to a
Fourier series for a closed curve, and in a similar way as aliasing
affects discrete Fourier series, samples of distinct functions may be
indistinguishable on a discrete grid when dealing with spherical
harmonics. High-frequency basis functions may alias to 
lower frequencies on a finite $p$-grid and introduce
pollution to lower frequencies which renders the simulation
unstable. This is a typical problem when dealing with non-linear
manipulation of functions (mean and Gaussian curvatures, surface
gradients, etc) and it is necessary to keep the aliasing
errors under control for these calculations. 

When computing geometrical quantities
involving non-linear manipulations, we will therefore first upsample the surface
and the functions involved to a finer grid\footnoteremember{upsampling_downsampling}{By upsampling we mean to pad with zeros in the spectral space and then transform back to the physical space, the finer grid. By downsampling back again we mean to transform to the spectral space, cut the extra modes and trasform back the physical space on the original grid.}, perform differentiation and
non-linear manipulations, and then downsample\footnoterecall{upsampling_downsampling} the result to the
original grid. The commonly used "$3/2$ upsampling rate" \cite{orszag} is not a good
choice for our purposes since highly non-linear functions are involved for which
this standard upsampling rate might not be sufficient. The adequate
upsampling rate depends on the specific shape of the surface, and we
will adopt an adaptive upsampling procedure introduced by Rahimian et
al. \cite{1} which is based on the decay of the mean curvature
spectrum.

\section{Numerical integration}
\label{num_int}

\subsection{Regular integrals and discrete orthogonality}
\label{regular_quadrature}
We will use a Gaussian spherical grid: in the longitudinal direction we discretize $\phi \in [0,2\pi)$ using $2p+2$ equidistant points $\left\{ \phi_i=\frac{\pi i}{p+1} \right \}_{i=0}^{2p+1}$, while in the non-periodic direction $\theta \in [0,\pi]$ we set $\left\{ \theta_j=\cos^{-1}(t_j)\right \}_{j=0}^{p}$, where $t_j \in [-1,1]$ are the Gauss-Legendre quadrature nodes with corresponding Gaussian quadrature weights $w_j^G$. This is a standard choice for spectrally accurate integration when using a spherical harmonics expansions \cite{1,6,20}.\\
Now we can define a quadrature rule for regular integrals:
\begin{equation}
\int_{\gamma} f d\gamma=\sum_{j=0}^p \sum_{k=0}^{2p+1} w_{j} f(\theta_j,\phi_k)W(\theta_j,\phi_k),
\end{equation}
where $w_{j}=\frac{\pi}{p}\frac{w_j^G}{ sin(\theta_j)}$ and $W(\theta_j,\phi_k)$ is the infinitesimal area element of the surface $\gamma$.\\
This choice of quadrature rule ensures the discrete orthogonality property: \\
\begin{equation}
\label{eq:discrete_ortho}
(Y_n^m,Y_{n'}^{m'})^p=\delta_{nn'}\delta_{mm'}, 
\end{equation}
for $0 \leq n,n' \leq p$ and $|m| \leq n, |m'| \leq n'$, see \cite{19} for details.\\
In this setting, the discrete version of eq. (\ref{eq:back_transf}) is:
\begin{equation}
\label{eq:back_transf_discr}
f_n^m=\int_0^\pi \int_0^{2\pi}f(\theta,\phi)\overline{Y_n^m(\theta,\phi)} \sin \theta d\theta d\phi \approx  \frac{\pi}{p} \sum_{i=1}^{(p+1)} \sum_{j=1}^{2(p+1)} w_j^G f(\theta_i,\phi_j)\overline{Y_n^m(\theta_i,\phi_j)} .
\end{equation}
We will refer to eq. (\ref{eq:back_transf_discr}) as regular quadrature.
Fast algorithms for computing backward and forward spherical harmonics transforms based on this grid are available \cite{fast_spharm1},\cite{7}.

\subsection{Singular integration}
\label{singular_integration}
When computing the integrals in (\ref{eq:main_eq}), for
$\textbf{x}_0=\textbf{x}$ we have a singularity. In order to compute the integrals we will make use of the following theorem \cite{8,9}:
\begin{theorem}
\label{teo_rot}
For any smooth function $f$ defined on a $C^\infty$ surface $\gamma$
globally parametrized by spherical coordinates $(\theta,\phi)$, the
harmonic potential evaluated at the north pole $\textbf{x}(0,0)$ is given by
\[
I[f](0,0)=\int_0^\pi \int_0^{2\pi}
\frac{f(\textbf{x}(\theta,\phi))W(\theta,\phi)}{|\textbf{x}(0,0)-\textbf{x}(\theta,\phi)|}
d\theta d\phi. 
\]
where $W$ represents the infinitesimal area element. 
For the quadrature rule 
\begin{equation}
\label{eq:quad_rule_sing}
Q[f](0,0)=  \sum_{j=0}^p \sum_{k=0}^{2p+1} \frac{w_j^s}{|\textbf{x}(0,0)-\textbf{x}(\theta_j,\phi_k)|}f(\textbf{x}(\theta_j,\phi_k))\frac{W(\theta_j,\phi_k)}{sin(\theta_j)} 
\end{equation}
where
\begin{equation}
w_j^s= w_j\sum_{n=0}^p 2\sin(\theta_j/2)P_n(\cos \theta_j),
\end{equation}
it holds that $Q$ converges superalgebraically to $I$ with $p$.
\end{theorem}
The theorem is based on the fact that the spherical harmonics are
eigenfunctions of the Laplace operator on the sphere. It can be
extended to the single and double layer kernels using symmetries and
smoothness properties as shown in \cite{6,1}. 
Using the quadrature rule (\ref{eq:quad_rule_sing}) we
can hence achieve an exponential convergence for the integrals that need to
be computed in (\ref{eq:main_eq}).

The theorem is stated for the target point $\textbf{x}_0$ at the north
pole of the parametrization, and does not hold in general. Hence, for
any target point, the coordinate system will be rotated such that this
quadrature rule can be applied. Analytical expressions for rotating the coefficients of the spherical harmonic expansion are given in
\ref{appendix2}. A fast algorithm for spherical grid
rotation with application to singular quadrature is given by Gimbutas and Veerapaneni \cite{2} where the computational complexity for all targets on the surface is $O(p^4\log{p})$.

\subsection{Nearly-singular integration}
\label{nearly-sing}

As two drops get closer, the integrals in Eq.(\ref{eq:main_eq}) get
more difficult to evaluate accurately due to the near singularities in
the integrals, see Fig. \ref{fig:NS_error_male1}. 
For this situation we have no analytical treatment as in the singular
case. The nearly singular integration in 3D remains an area of active
research, both for fixed and time-dependent geometries. In the first case, data can be precomputed and stored in order to accelerate the algorithm, but this is not possible in our case when the geometry is moving and deforming. Most of the methods available are based on
refinement and/or interpolation techniques: Zhao
et. al. \cite{21} proposed a technique based on the method of floating
partition of unity for both singular and nearly singular integrals; Bremer and Gimbutas \cite{Bremer2012,Bremer2013} describe a modified Nystr\"om method for the discretization of weakly singular integral operators by first applying a suitable transformation to simplify the integrand and then using a precomputed table of quadrature rules to evaluate the resulting integral efficiently;
Tlupova and Beale \cite{20}
developed a third order method using a regularized kernel, and, more recently, af Klinteberg and Tornberg
\cite{Ludvig_QBX} presented a fast and accurate method based
on QBX (quadrature-by-expansion) for spheroidal surfaces, utilizing precomputation for efficiency. Ying, Biros and Zoris \cite{13} proposed a procedure based on
interpolation to compute the nearly-singular integral.
We will adapt and refine the last approach for our application.

When the drops are well separated, at some distance larger than $\tilde{A}(h)$ from each other, the regular
quadrature works well. As they get closer, regular
quadrature on a finer grid can still be used. Here, the density is interpolated to the finer (upsampled) grid, where the nearly singular kernel is better resolved. But at some point, a special quadrature
method is needed since the quadrature error grows exponentially as we
approach the surface \cite{Ludvig_QBX} and it is not possible to resolve
the problem by grid refinement, i.e. upsampling. 

To formalize this division into regions, let $\bold{x}_0$ denote the
target point and $h=\frac{R\pi}{p+1}$ the grid size, where $R$ is the mean radius of the drop. 
We consider the following three regions:
\begin{itemize}
\item Well separated: $\Omega_0=\{\bold{x}_0 \in \gamma_i \ | \ dist(\textbf{x}_0,\gamma^* \setminus \gamma_i)\in(\tilde{A}(h),\infty), \ i=1,\dots,N \}$:\\
A regular quadrature rule on the standard grid is used;
\item Intermediate: $\Omega_1=\{\bold{x}_0 \in \gamma_i \ | \ dist(\textbf{x}_0,\gamma^* \setminus \gamma_i)\in(h,\tilde{A}(h)], \ i=1,\dots,N \}$:\\
A regular quadrature rule on the grid upsampled by a factor $\tilde{U}$ is used;
\item Nearest: $\Omega_2=\{\bold{x}_0 \in \gamma_i \ | \ dist(\textbf{x}_0,\gamma^* \setminus \gamma_i)\in(0,h], \ i=1,\dots,N \}$\\
In this case using the upsampled grid is not enough and we need a
special technique to compute the nearly-singular integral, as will be
described below. 
\end{itemize}

\begin{figure}[htbp]
	\centering
		\includegraphics[width=100mm]{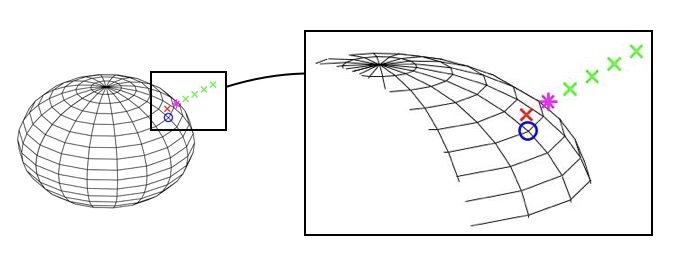}
		
\caption{Schematic picture of the special quadrature for the nearly singular integrals. The pink star represents the target point $\textbf{x}_0$, the blue circle is the closest grid point used as starting point for the Newton iteration (\ref{eq:newton_ite}). The
  red cross represents $\textbf{x}_*$, the projection on the surface
  of the  target point $\textbf{x}_0$, and the green crosses are the points far away from the surface used for the 1D interpolation.}
\label{fig:sph_ns23}
\end{figure}

The idea introduced in Ying et al. \cite{13} for the nearly singular
integration was to find the point $\textbf{x}_*$ on the surface that is closest to
the target point $\textbf{x}_0$. Then, continuing along a line that
passes through $\textbf{x}_*$ and $\textbf{x}_0$, the integral is
evaluated at a number of points $\textbf{x}_1,\dots,\textbf{x}_n$ further away from the surface (see Figure \ref{fig:sph_ns23}). This can be done by regular quadrature on the standard grid or on the upsampled grid, depending if the point belongs to $\Omega_0$ or $\Omega_1$. The value of the integral on the surface (at $\textbf{x}_*$) needs to be computed by special quadrature for singular integrals. At this point a 1D Lagrangian interpolation is used to compute the value at $\textbf{x}_0$ by interpolating the values at $\textbf{x}_*$ and $x_i, \ i=1,\dots,N$. 

We will follow the same approach as above in the following way. A cell list algorithm has been implemented to find in which region the target point $\textbf{x}_0$ belongs. For the regions $\Omega_0$ and $\Omega_1$, the regular quadrature rule has already been introduced in Section \ref{regular_quadrature}. For the target points belonging to the region $\Omega_2$, we need to find the closest point on the surface, $\textbf{x}_*$. Using the cell list mentioned above, we can also find the closest grid point to $\textbf{x}_0$ (see Fig. \ref{fig:sph_ns23}); then we use it as starting point for a Newton iteration to find the closest point on the surface (not necessarily a grid point). We need to minimize the distance between the closest drop surface and the target point. Since the surface is represented by a spherical harmonic expansion, we want to find the parameters $(\phi,\theta)$ corresponding to the point which minimizes the distance
\[\min_{(\phi,\theta)}\Arrowvert \bold{x}(\phi,\theta)-\bold{x}_0 \Arrowvert=\min_{(\phi,\theta)}\Arrowvert \sum_{n=0}^p \sum_{m=-n}^n \bold{x}_n^m Y_n^m(\phi,\theta)-\bold{x}_0 \Arrowvert\]
(where $\bold{x}_n^m$ are the spherical harmonic coefficients of the drop surface, eq. (\ref{eq:surface_ev}).)\\
To solve this optimization problem, we look for the zeros of the
gradient of $$f(\phi,\theta)=\Arrowvert
\bold{x}(\phi,\theta)-\bold{x}_0 \Arrowvert,$$ where we have
\begin{equation}
\label{eq:f_phi}
\frac{\partial f}{\partial \phi}=2(\sum_n \sum_m \bold{x}_n^m Y_n^m-\bold{x}_0)\cdot(\sum_n \sum_m \bold{x}_n^m \frac{\partial Y_n^m}{\partial \phi})=0
\end{equation}
\begin{equation}
\label{eq:f_theta}
\frac{\partial f}{\partial \theta}=2(\sum_n \sum_m \bold{x}_n^m Y_n^m-\bold{x}_0)\cdot(\sum_n \sum_m \bold{x}_n^m \frac{\partial Y_n^m}{\partial \theta})=0
\end{equation}
We solve the system (\ref{eq:f_phi})-(\ref{eq:f_theta}) using a Newton iterative method,
\begin{equation}
\label{eq:newton_ite}
(\phi^{n+1},\theta^{n+1})=(\phi^{n},\theta^{n})-H^{-1}\nabla(f(\phi^{n},\theta^{n}))
\end{equation}
where H is the Hessian of $f$.\footnote{Note that the inverse of H
  might be not well defined. To avoid this problem, an approach
  exploited in the Levenberg-Marquardt algorithm (which uses an
  approximate Hessian) is to add a scaled identity matrix to the
  Hessian, $\mu I$. This results in slower but more reliable
  convergence where the Hessian doesn't provide useful information
  \cite{16}.}

Once $\textbf{x}_*$ has been determined, we could apply Theorem
\ref{teo_rot} for computing the integral at this location.  We have
however already computed the single and double layer integrals for the
\textsl{self-self} interaction with all grid points on the surface as target points,
and these discrete values in term determine a spherical harmonics
expansion of the field over the surface. Hence, we can simply use
this expansion to determine the value at $\textbf{x}_*$, and this same
expansion can be used for any target point close to the surface where the
special quadrature is invoked. This approach is computationally cheaper than the direct application of Theorem \ref{teo_rot} even if it is used only for a single point.
Finally, the interpolation procedure previously described to obtain the
potential at the target point $\textbf{x}_0$ can be used.

In this procedure, there are a number of parameters to set,
\begin{itemize}
\item $\tilde{A}$, the maximum distance before we need to upsample the grid for the regular quadrature,
\item $\tilde{U}$, the upsampling rate used in the intermediate
  region,
\item $L$, the number of points used for interpolation for target
  points in the nearest region,
\item $\tilde{D}$, the distance between these points.
\end{itemize}
Clearly all these parameters depend on the resolution, so we will express them in terms of $h$, and we will see in the next section how they affect the accuracy.

\subsubsection{Testing the nearly singular quadrature}
\begin{figure}
    \centering
          \hspace{-4cm}
     \begin{subfigure}[b]{0.3\textwidth}
        \includegraphics[width=\textwidth]{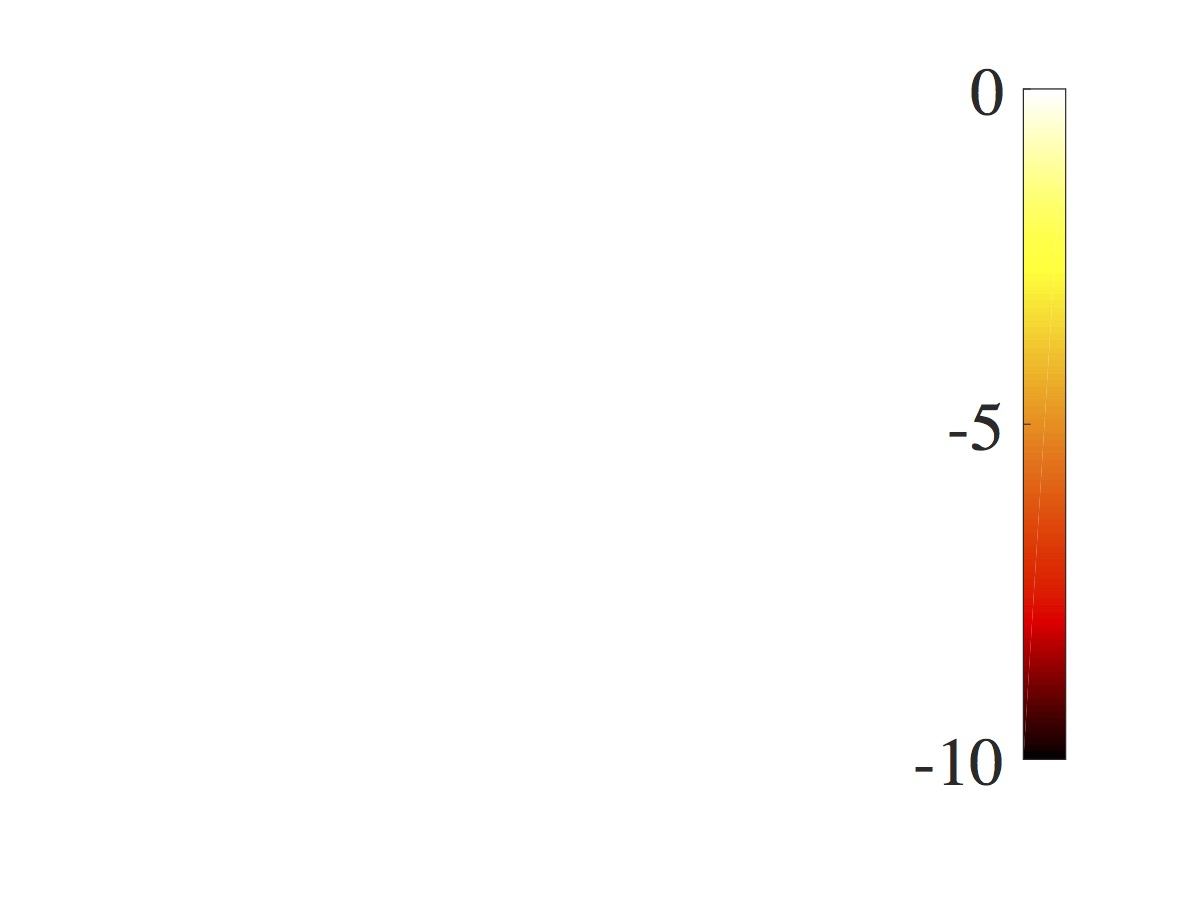}
    \end{subfigure}
     \hspace{-0.7cm}
    \begin{subfigure}[b]{0.38\textwidth}
        \includegraphics[width=\textwidth]{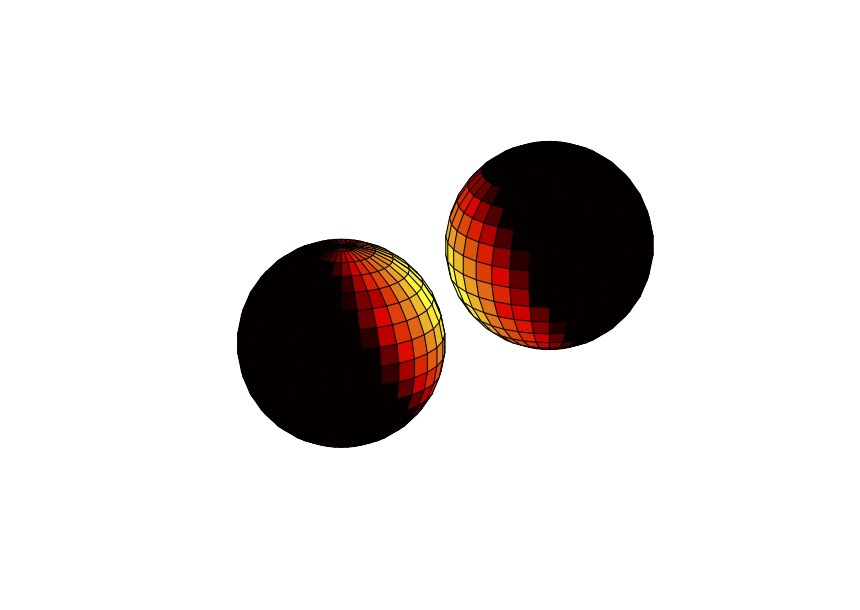}
        \caption{}
        \label{fig:NS_error_male1}
    \end{subfigure}
    ~ 
      \hspace{-2cm}
      \begin{subfigure}[b]{0.38\textwidth}
        \includegraphics[width=\textwidth]{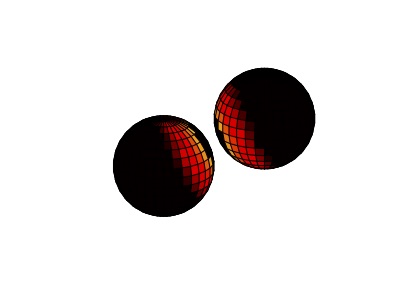}
        \caption{}
        \label{fig:NS_error_male2}
    \end{subfigure}
     \hspace{-1.6cm}
    \begin{subfigure}[b]{0.38\textwidth}
        \includegraphics[width=\textwidth]{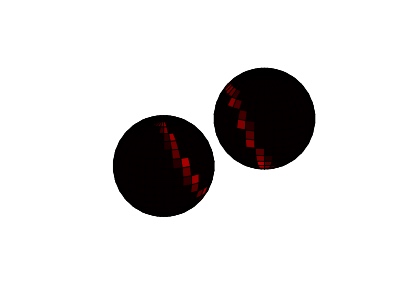}
        \caption{}
        \label{fig:NS_error_male3}
    \end{subfigure}
    \caption{Error (log scale) in computing the double-layer identity, $p=15$. (a) Regular quadrature; (b) Special quadrature with $\tilde{A}=\sqrt{h}, \tilde{U}=4$; (c) Special quadrature with $\tilde{A}=5h, \tilde{U}=4$.}\label{fig:NS_error_male}
\end{figure}

We will first test the nearly-singular quadrature algorithm for the double-layer potential on a unit sphere; in Fig. \ref{fig:NS_error_male} it is shown how the error behaves with and without special quadrature. We notice that the upsampled region is not big enough if we fix the distance $\tilde{A}=\sqrt{h}$  as suggested in \cite{13} (Fig. \ref{fig:NS_error_male2}), but $5h$ seems to be safer to ensure good accuracy (Fig. \ref{fig:NS_error_male3}). For this specific case the error is computed for the double-layer identity, where we know the analytical solution:
\begin{equation*}
\int_\Gamma T_{ijk}(x,y)\hat{n_k}(y)ds(y)=\delta_{ij} 
\begin{cases} 
0 &\mbox{if } x $ outside $ \Gamma, \\ 
4\pi &\mbox{if } x \in \Gamma, \\ 
8\pi &\mbox{if } x $ inside $ \Gamma. \\ 
\end{cases}
\end{equation*}

\begin{figure}
    \centering
         \hspace{-0.8cm}
    \begin{subfigure}[b]{0.45\textwidth}
        \includegraphics[width=75mm]{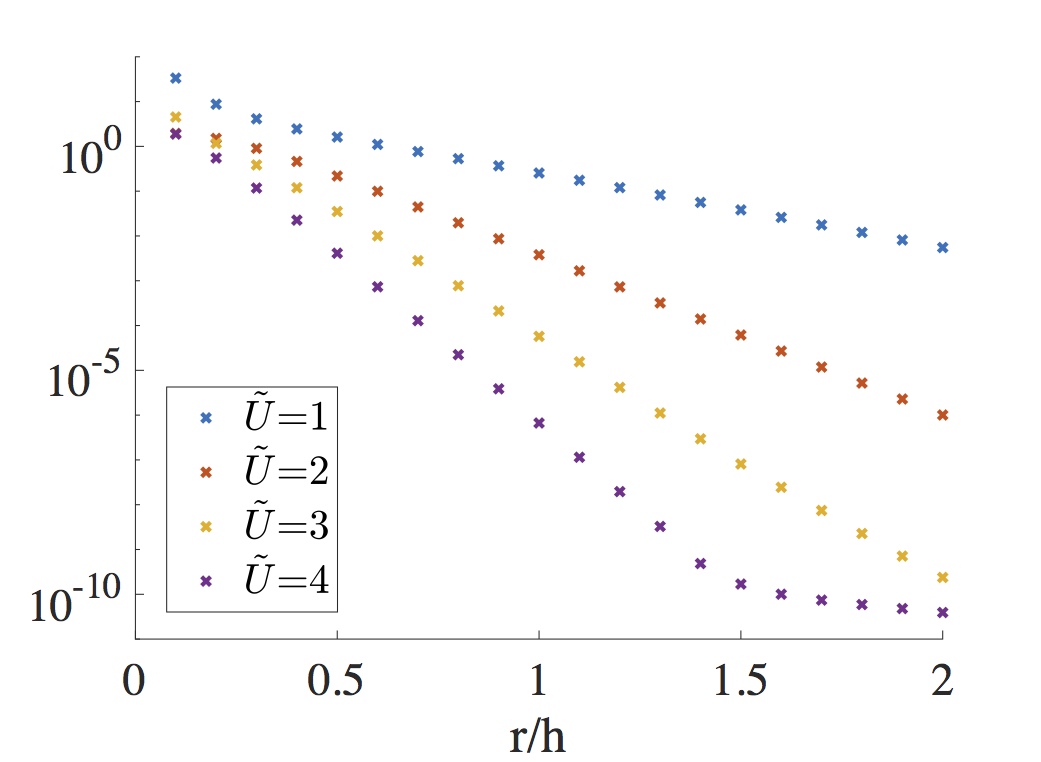}
        \caption{}
        \label{fig:jeffery1}
    \end{subfigure}
    ~ 
      \hspace{0.8cm}
      \begin{subfigure}[b]{0.45\textwidth}
        \includegraphics[width=75mm]{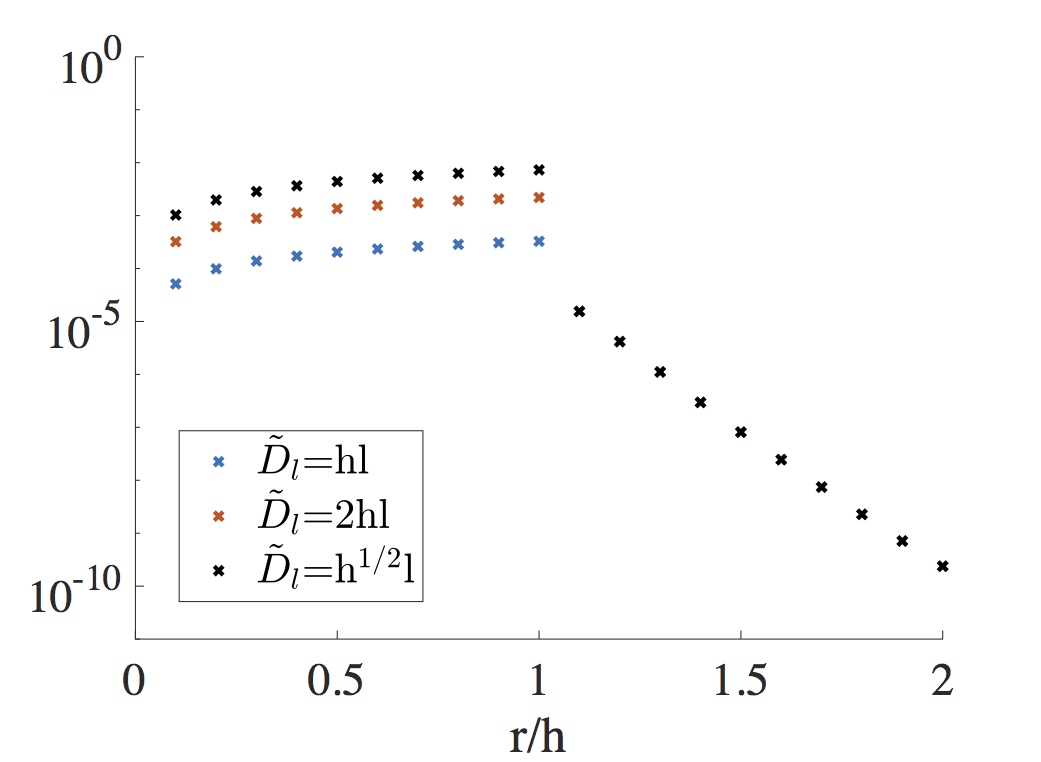}
        \caption{}
        \label{fig:jeffery2}
    \end{subfigure}
    
   \hspace{-0.8cm}
    \begin{subfigure}[b]{0.45\textwidth}
        \includegraphics[width=75mm]{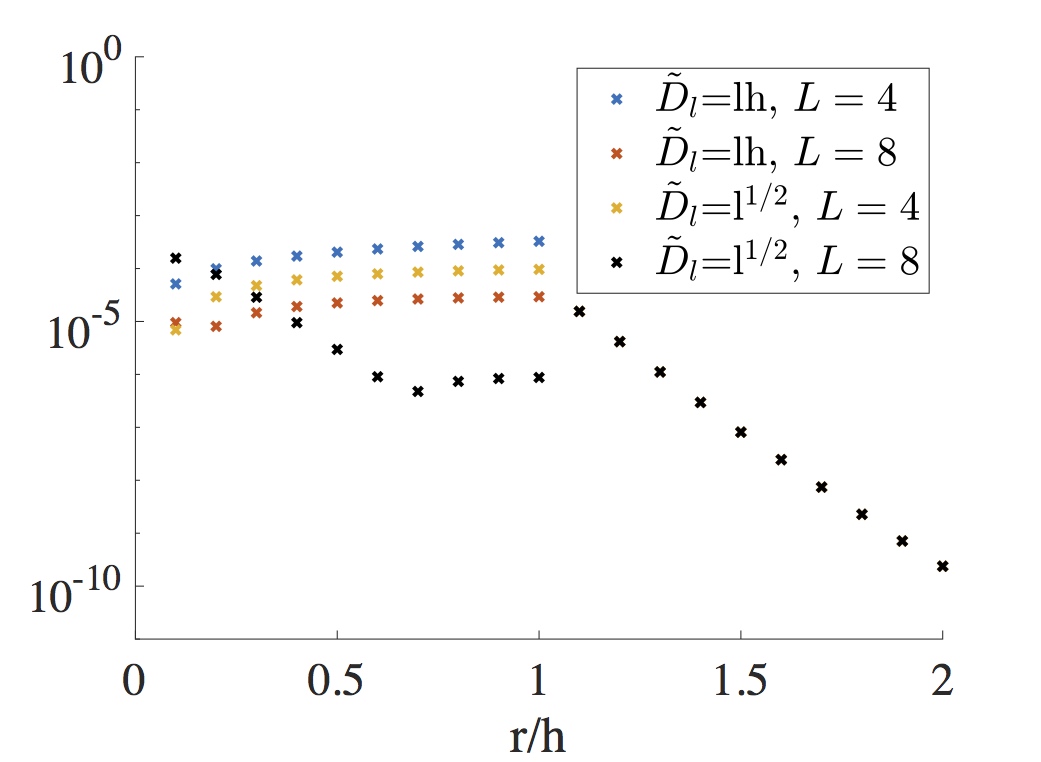}
         \caption{}
        \label{fig:jeffery3}
    \end{subfigure}
    ~ 
      \qquad
      \begin{subfigure}[b]{0.45\textwidth}
        \includegraphics[width=75mm]{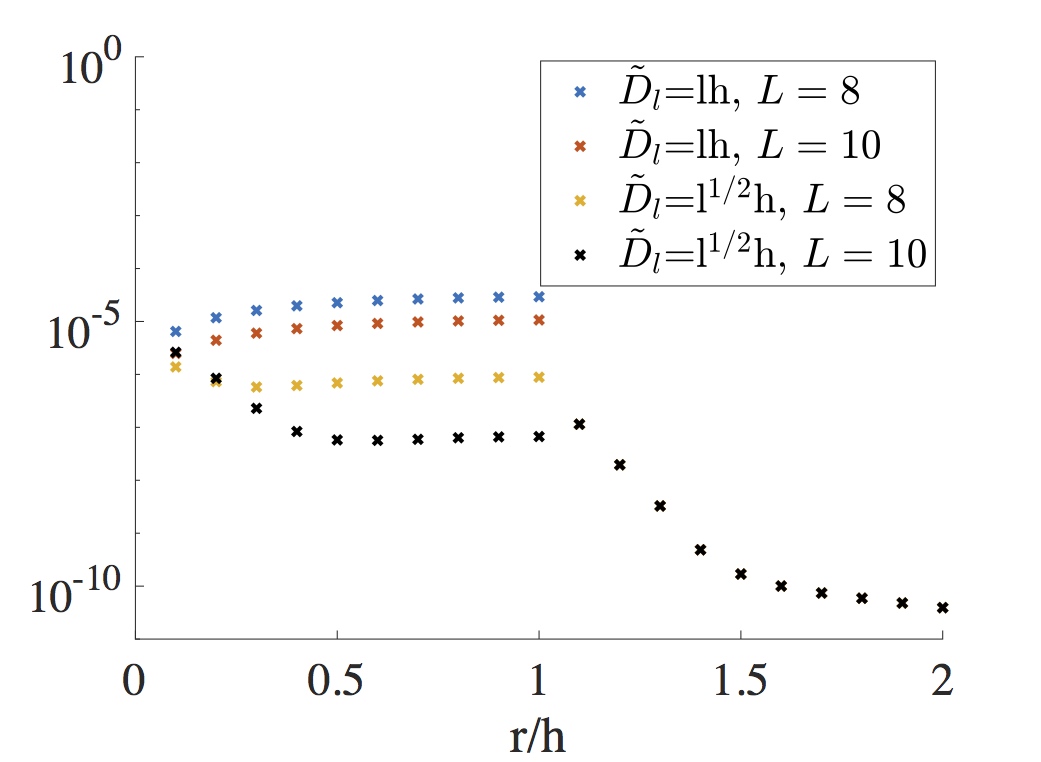}
        \caption{}
        \label{fig:jeffery4}
    \end{subfigure}
 
  \caption{$L_\infty$ norm of the error for the analytical test
      given by \cite{14}-\cite{15} vs the relative distance $r/h$;
      $r$ is the distance between the target point and the ellipsoid
      and $h=\frac{\pi}{p+1}$. $p=39$. (a) Regular quadrature with different upsampling rates; (b) Special quadrature with $L=4$ uniformly distributed interpolating points in the \textit{Nearest region} and regular quadrature on the upsampled grid with $\tilde{U}=3$ in the \textit{Intermediate region}; (c) Special quadrature with $L=4,8$ uniformly and non-uniformly distributed interpolating points in the \textit{Nearest region} and regular quadrature on the upsampled grid with $\tilde{U}=3$ in the \textit{Intermediate region}; (d) Special quadrature with $L=8,10$ uniformly and non-uniformly distributed interpolating points in the \textit{Nearest region} and regular quadrature on the upsampled grid with $\tilde{U}=4$ in the \textit{Intermediate region}.}
    \label{fig:Jeffery1234}
\end{figure}

This solution is constant outside the surface, and therefore
inadequate for studying the convergence of the interpolation method
(sketch in Fig. \ref{fig:sph_ns23}), so we will use instead an
analytical non-constant solution given by Jeffery \cite{14,15}
for a spheroid with aspect ratio (1:2) to analyze the method. 
In Fig. \ref{fig:Jeffery1234} we compare the error obtained by regular
quadrature with the one obtained by the special quadrature tuning differently the parameters $\tilde{U}$, $\tilde{D}_l$ and $L$. In Fig. \ref{fig:jeffery2} we test uniformly distributed points for the interpolation, where the distance between points is constant and depends only on the
resolution parameter $h$; in this sense the distribution of points is
linear in $l$. We tried different choices: $\tilde{D}_l=hl$,
$\tilde{D}_l=2hl$, $\tilde{D}_l=\sqrt{h}l$, $l=1,\dots,L$, but 
in none
of these cases the accuracy is sufficient, indeed we observe
a jump between the regular quadrature with upsampling and the special
quadrature. We observe that using a higher upsampling rate $\tilde{U}$ always gives lower errors. Also when we consider a sequence of interpolation points that scales with the distance as $\tilde{D}_l(h)=h\sqrt{l}$ and increase the number of interpolation points, $L=8,10$, we can achieve higher accuracy but we have to pay attention in tuning these parameters correctly,  as shown in Fig. \ref{fig:jeffery3}-\ref{fig:jeffery4}. 

\section{Time-stepping}
\label{time}
The choice of the time-stepping scheme for updating the drop positions and the surfactant evolution has to consider two main aspects: the computational cost and the control of error. The main numerical cost comes from solving the integral equation to obtain the velocity of the drops as presented in the previous Sections. We want to maintain the number of Stokes evaluations at a minimum; at the same time the numerical scheme for the time evolution needs to control the error in both quantities (drop/surfactant) and the time-step size should be adaptive balancing this two factors. \\ 
We recall that the system of ODEs (drop/surfactant) we are considering is given by
 \begin{equation}
 \label{eq:system1fefi}
\begin{cases}
\frac{d\textbf{x}}{dt}=\textbf{u}(\textbf{x},\sigma), \\
\frac{d\Gamma}{dt}=f_E(\textbf{x},\textbf{u},\Gamma)+f_I(\textbf{x},\Gamma), 
\end{cases}
 \end{equation}
where the velocity $\textbf{u}$ at the interface is given by the solution
of the boundary integral equations (\ref{eq:main_eq}), $f_E$ and $f_I$ represents respectively the convective and the diffusion component of the equation for the surfactant (see eq. (\ref{eq:surfactant0})),
$$f_E=-\nabla_\gamma \cdot (\Gamma \textbf{u}_\gamma) -2K(\textbf{x})\Gamma (\textbf{u}\cdot \textbf{n}),$$
$$f_I=\frac{1}{Pe}\nabla_\gamma^2\Gamma.$$
Due to the different nature of these two terms, we would like to treat them differently: an implicit scheme is a good choice for the diffusion term and an explicit method for the convective term \cite{IMEX}.

We will first introduce three explicit schemes for the drop evolution (Section \ref{drop_evolution}), then an IMEX scheme for the surfactant concentration (Section \ref{surfactant_evo}), with particular attention to the treatment of the implicit term; the coupling drop/surfactant is briefly explained in Section \ref{coupling_comp} where we also compare the different schemes considered.

\subsection{Drop evolution}
\label{drop_evolution}
For the drop evolution we compare three different adaptive Runge-Kutta schemes, with corresponding Butcher tableaus given by: 
\begin{center}
\hspace{-5em}
\begin{minipage}[t]{\widthof{\textbf{Explicit midpoint method} ($2^{nd}$ order)}}
\centering
\textbf{Explicit midpoint}\\($2^{nd}$ order)\\[2ex]
$\begin{array}{c|cc}
0 \\
1/2 & 1/2 \\
\hline 
&0&1\\
&1&0
\end{array}$
\end{minipage}\hspace{-5em}
\begin{minipage}[t]{\widthof{\textbf{RK23} (Bogacki-Shampine, $3^{rd}$ order)}}
\centering
\textbf{RK23} \\(Bogacki-Shampine, $3^{rd}$ order)\\[2ex]
$\begin{array}{c|ccccc}
0 \\
1/2 & 1/2 \\
3/4 & 0 & 3/4 \\
1 & 2/9 & 1/3 & 4/9 \\
\hline 
 & 2/9 & 1/3 & 4/9 \\
 & 7/24 & 1/4 & 1/3 & 1/8 \\
\end{array}$
\end{minipage}\hspace{-2em}
\begin{minipage}[t]{\widthof{\textbf{Kutta's method} ($3^{rd}$ order)}}
\centering
\textbf{Kutta's method} \\($3^{rd}$ order)\\[2ex]
$\begin{array}{c|ccccc}
0 \\
1/2 & 1/2 \\
1 & -1 & 2 \\
\hline 
 & 1/6 & 4/6 & 1/6 \\
 & 0 & 1  \\
\end{array}$
\end{minipage}
\end{center}
where both the used scheme (higher order) and the scheme of lower order used for computing the local error are shown.

\subsection{Surfactant evolution}
\label{surfactant_evo}
For the surfactant evolution we consider a $2^{nd}$ order IMEX Runge-Kutta Scheme \cite{IMEX} with time step size $dt$ as follows:

\begin{itemize}
\item $\tilde{k}_0=f_E(\textbf{x}^{(t)},\textbf{u}^{(t)},\Gamma^{(t)})$;
\item solve for $\Gamma^{(t+dt/2)}$ (see Section \ref{preconditioning} for more details):\\
\begin{equation}
\label{eq:systemforgammatpdt2}
\Gamma^{(t+dt/2)}=\Gamma^{(t)}+dt/2[f_E(\textbf{x}^{(t)},\textbf{u}^{(t)},\Gamma^{(t)})+\frac{1}{Pe}\nabla_S^2 \Gamma^{(t+dt/2)}],
\end{equation}
and compute\\
$k_{1/2}=f_I(\textbf{x}^{(t+dt/2)},\Gamma^{(t+dt/2)})$;
\item $\tilde{k}_{1/2}=f_E(\textbf{x}^{(t+dt/2)},\textbf{u}^{(t+dt/2)},\Gamma^{(t+dt/2)})$;
\end{itemize}
\begin{equation}
\label{eq:gammatpiudt}
\Gamma^{(t+dt)}=\Gamma^{(t)}+dt(\tilde{k}_{1/2}+k_{1/2}).
\end{equation}
The input data $\textbf{u}^{(i)}$, for $i=t,t+dt/2$, and the corresponding position vectors $\textbf{x}^{(i)}$,  are given by the time-solver for the drop evolution (see Algorithm \ref{alg:pseudo_code}).\\
Since we want an adaptive method able to keep the error below a fixed tolerance, we compare two methods to control the surfactant evolution:
\begin{enumerate}
\item Error computed with IMEX1 ($1^{st}$ order):
\begin{equation}
\Gamma^{(t+dt)}_{IMEX1}=\Gamma^{(t)}+dt(f_E(\textbf{x}^{(t)},\textbf{u}^{t},\Gamma^{(t)})+f_I(\textbf{x}^{(t+1)},\Gamma^{(t+dt)}_{IMEX1}));
\end{equation}
in this case the estimation of the error is given by
\begin{equation}
\label{eq:err_IMEX1}
err_{surfactant}^{(t+dt)}=\frac{\| \Gamma^{(t+dt)}-\Gamma^{(t+dt)}_{IMEX1}\|_{\infty}}{\| \Gamma^{(t+dt)}\|_{\infty}}.
\end{equation}
\item Error based on the local conservation of surfactant (\textit{conservation error}): in the case of insoluble surfactant the mass is conserved on the interface,
$$\frac{d}{dt} \int_{\gamma(t)}\Gamma dS=0.$$
As an estimate for the surfactant error, we compute the local error in conservation of mass
\begin{equation}
\label{eq:err_conservation}
err_{surfactant}^{(t+dt)}=\frac{|\int_{\gamma(t+dt)}\Gamma^{(t+dt)}dS-\int_{\gamma(t)}\Gamma^{(t)}dS|}{|\int_{\gamma(t)}\Gamma^{(t)}dS|},
\end{equation}
where by \textsl{local} we refer to the change in mass between the current and the previous time-step and where the integral is evaluated numerically as explained in Section \ref{regular_quadrature}.
\end{enumerate}
We also considered a Semi-implicit Spectral Deferred Correction Method \cite{SDC1,SDC2} but it turned out to be too expensive for our problem since it requires several BIM solver calls to achieve the desired tolerance for the coupled system.

\subsubsection{Preconditioning}
\label{preconditioning}
Using the IMEX scheme for the surfactant, we need to solve the system generated by eq. (\ref{eq:systemforgammatpdt2}) for every time step at the stage $(t+dt/2)$,
which we can write in a compact form as:
\begin{equation}
A\textbf{y}=\textbf{b}
\end{equation}
where
$$A=I-\frac{dt}{2}\frac{1}{Pe}\nabla_S^2,$$
$$\textbf{b}=\Gamma^{(t)}+\frac{dt}{2}f_E(\textbf{x}^{(t)},\textbf{u}^{(t)},\Gamma^{(t)}),$$
and $\textbf{y}=\Gamma^{(t+dt/2)}$. In practice we actually solve the system for 
the spherical harmonics coefficients, see eq. (\ref{eq:galerkin_surfactant}).

\begin{figure}
{\includegraphics[width=6cm]{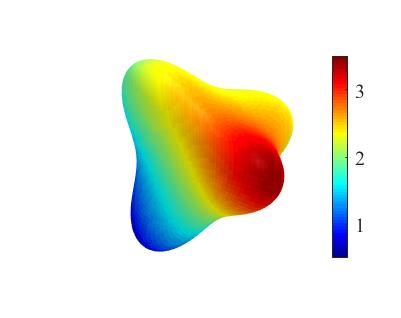}}
\caption{Surface used to test the preconditioner; the colors denote the surfactant concentration.}
\label{fig:weird_surface}
\end{figure}

The conditioning of this system gets worse with increased diffusion,
so we need a preconditioner to accelerate the convergence. It is well
known that the spherical harmonics are eigenfunctions of the
Laplace-Beltrami operator on the sphere \cite{SPHARM},
\begin{equation}
\Delta_{\mathbb{S}^2} Y_n^m=-n(n+1)Y_n^m.
\end{equation}
We use the diagonal matrix resulting on a sphere as a preconditioner for a general problem (general surface). In the following plots we are comparing the number of iterations (gmres tol=10e-12) for different diffusion coefficients (different $Pe$) with initial concentration $\Gamma=2+\textbf{x}$, on the surface\footnote{The surface is given by:
\begin{equation*}
\begin{cases}
		x_1(\theta,\phi)= \rho(\theta,\phi) cos(\phi) sin(\theta);\\
        x_2(\theta,\phi)= \rho(\theta,\phi) sin(\theta) sin(\phi);\\
        x_3(\theta,\phi)= \rho(\theta,\phi) cos(\theta)
\end{cases}
\end{equation*}
with $\rho(\theta,\phi)=0.7+0.3e^{-3Re(Y_3^2(\theta,\phi))}$.
} showed in Fig. \ref{fig:weird_surface}.

\begin{figure}[htbp]
\begin{subfigure}[b]{0.3\textwidth}
	\centering		
		\includegraphics[width=70mm]
{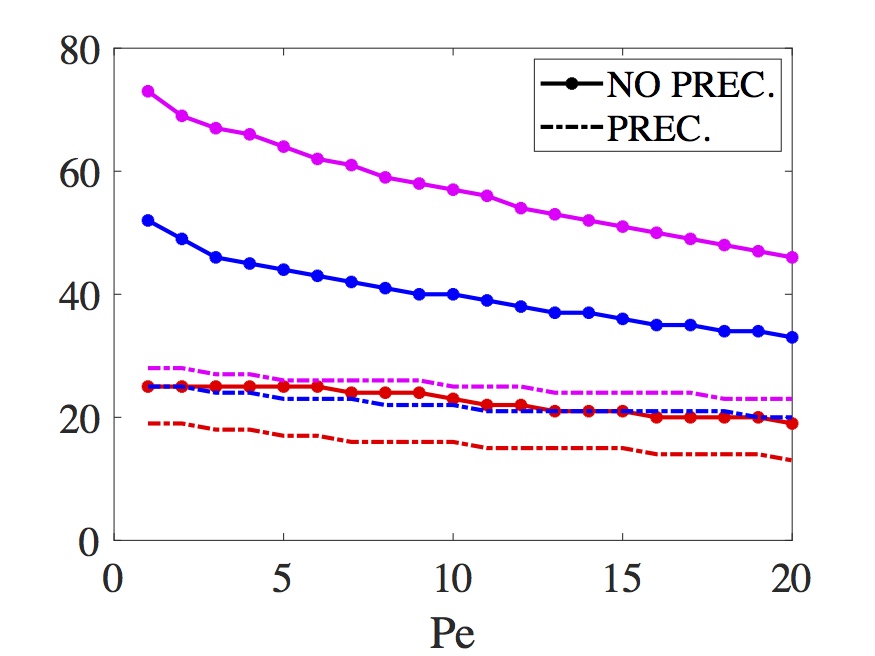}
\caption{Number of iterations for the IMEX scheme}
\label{fig:n_ite}
\end{subfigure}\hspace{2cm}
\begin{subfigure}[b]{0.4\textwidth}
	\includegraphics[width=70mm]{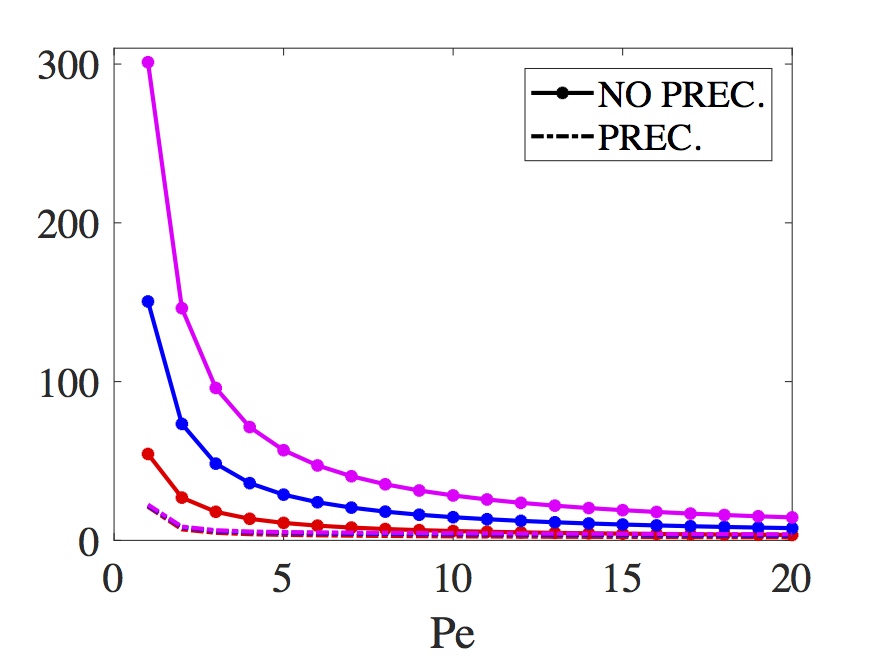}
\caption{Condition number for the surfactant matrix of the IMEX scheme}
\label{fig:cond_num}
\end{subfigure}
\caption{The effect of the preconditioner for the IMEX scheme for different resolutions: $p=9,15,21$ respectively represented by the red, blue and magenta colors. When applying the preconditioner introduced in Section \ref{preconditioning} the number of iterations needed is strongly reduced as well as the condition number, which becomes of the same order independent of $p$.}
\end{figure}

In Fig. \ref{fig:n_ite} we show a substantial reduction in the number
of iterations needed for the preconditioned algorithm and a gain increasing with the order of spherical harmonics; we also
show (Fig. \ref{fig:cond_num}) the condition number of the matrix $A$
compared to the condition number of the preconditioned matrix
$B=S^{-1}A$, where S is the diagonal matrix obtained on the
sphere, and we observe that, when the preconditioner is used, the
order of the condition number remains of the same order for all the different orders of spherical harmonics expansion considered.

\subsection{Coupling and comparison}
\label{coupling_comp}
The explicit scheme for the drop evolution and the IMEX scheme for the surfactant concentration evolution need to be coupled together, as shown in the pseudo-code Algorithm \ref{alg:pseudo_code}. We would like the drop surface and the surfactant concentration evolutions to use the same temporal stages, and this is the basic reason why we considered and restricted the comparison to the above mentioned schemes. Indeed, they all use the same stages with an exception given by the case of RK23 where we would need an extra computation of the surfactant at $t+\frac{3}{4}dt$. This can be done in two ways:
either we half the time-step for the surfactant (but we would do an extra unnecessary computation in $t+dt/4$), or we simply compute the solution in $t+\frac{3}{4}dt$ using a first order order IMEX method with time-step $dt/4$, i.e. we solve for $\Gamma^{t+\frac{3}{4}dt}$:
$$\Gamma^{(t+\frac{3}{4}dt)}=\Gamma^{(t+dt/2)}+\frac{dt}{4}(f_E(\textbf{x}^{(t+dt/2)},\textbf{u}^{(t+t/2)},\Gamma^{(t+dt/2)})+f_I(\textbf{x}^{(t+\frac{3}{4}dt)},\Gamma^{(t+\frac{3}{4}dt)})).$$
We adopt the second choice.

\begin{figure}[htbp]
	\centering
	  \begin{subfigure}[b]{0.45\textwidth}
		\includegraphics[width=80mm]{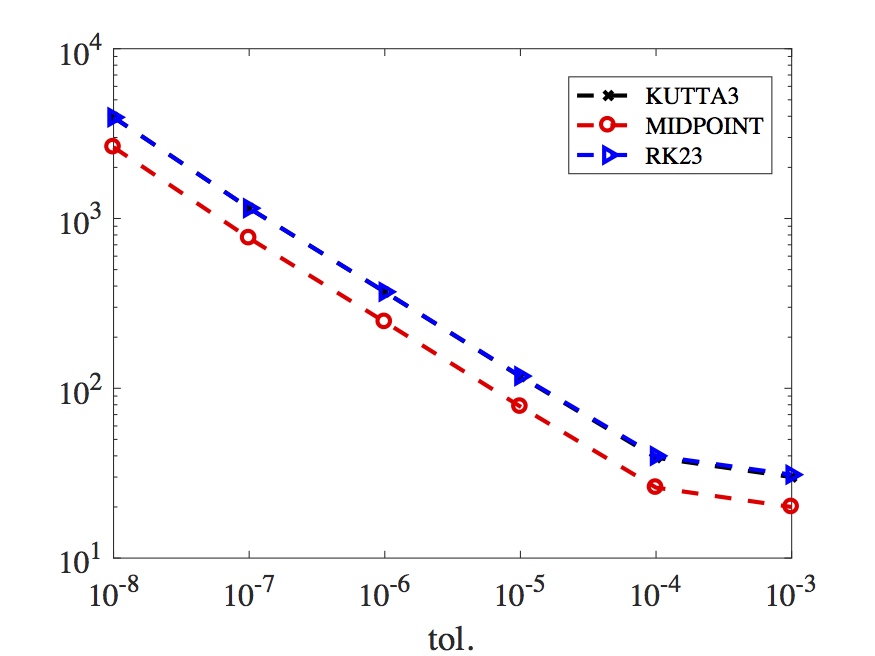}
		\caption{IMEX1, eq. (\ref{eq:err_IMEX1})}
		 \label{fig:tol_stokes_IMEX1}
		\end{subfigure}
		\quad
		 \begin{subfigure}[b]{0.45\textwidth}
		\includegraphics[width=80mm]{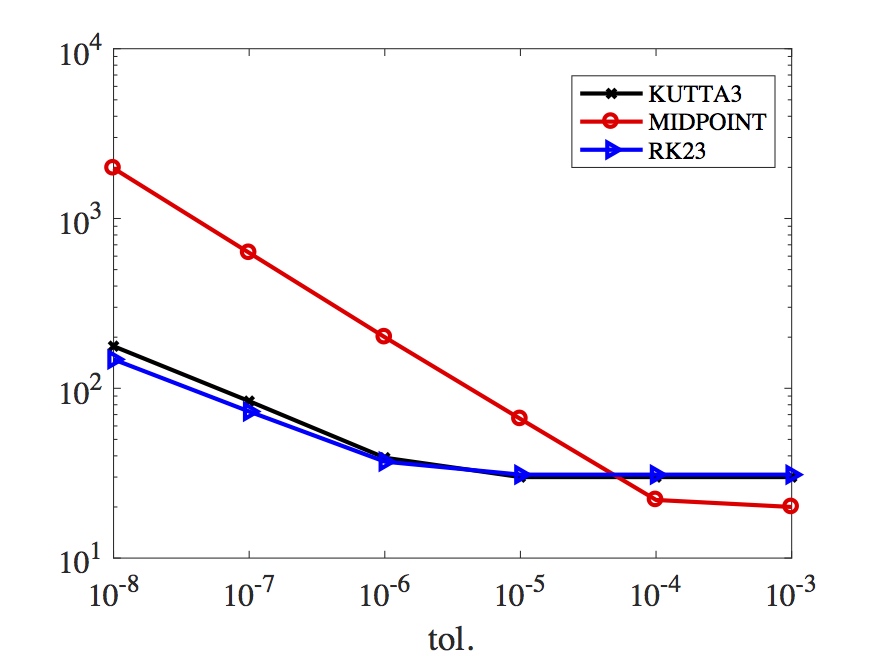}
     \caption{Conservation error, eq. (\ref{eq:err_conservation})}
      \label{fig:tol_stokes_cons}
		\end{subfigure}
		\caption{The number of Stokes evaluations plotted versus tolerance using the two different estimates for the surfactant error. For this simulation we used: $p=21$, final time $T=1$, initial shape: ellipsoid (a/c=1/2) with initial surfactant concentration $\Gamma=1+x$ immersed in a shear flow with $Ca=0.2$, $\lambda=1$, $Pe=10$, $E=0.2$, $x_s=0.3$.}
  \label{fig:tol_stokes}
\end{figure}

\begin{figure}[htbp]
	\centering
	\begin{subfigure}[b]{0.45\textwidth}
		\includegraphics[width=80mm]{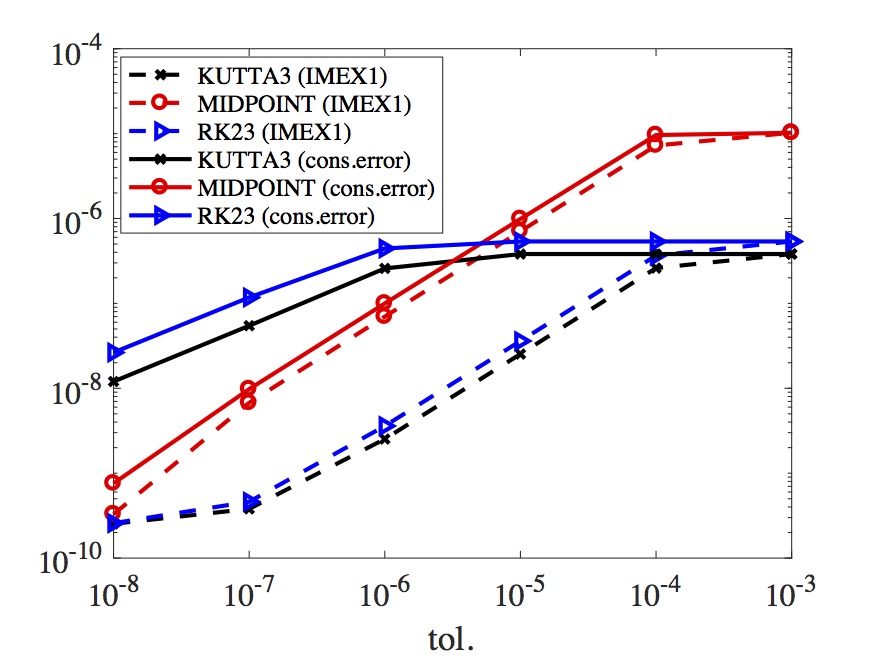}
		\caption{Drop}
		  \label{fig:tol_errors_drop}
		\end{subfigure}
		\quad
		\hspace{0.1cm}
		\begin{subfigure}[b]{0.45\textwidth}
		\includegraphics[width=80mm]
{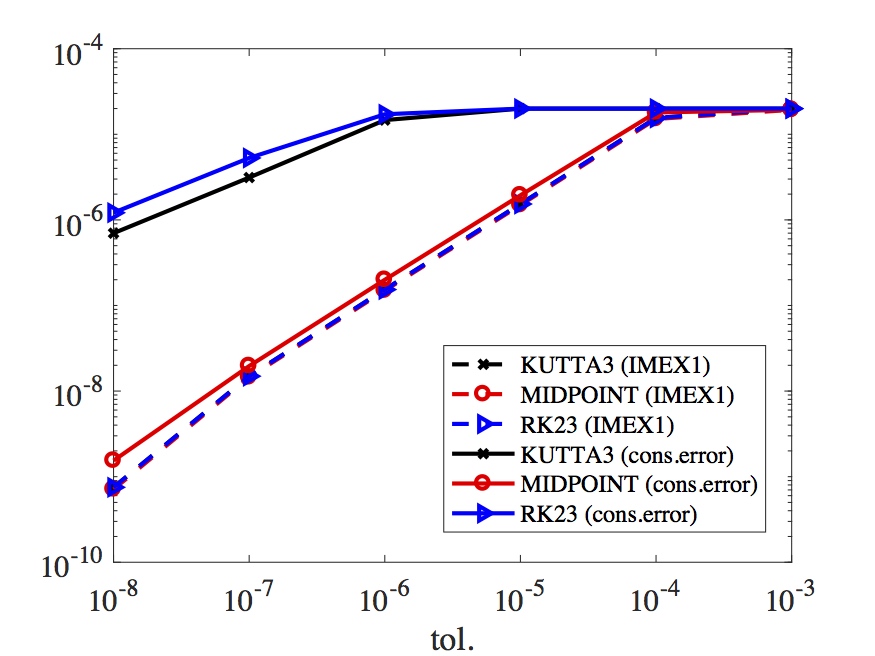}
\caption{Surfactant}
  \label{fig:tol_errors_surf}
		\end{subfigure}
\caption{Relative error vs tolerance. We compare the error (L2-norm) for the drop position and surfactant concentration. The data of this simulation are the same as in Fig. \ref{fig:tol_stokes}, the reference solution is computed using the Midpoint scheme combined with the second order IMEX scheme and a fixed time-step $\Delta t=10^{-5}$.}
  \label{fig:tol_errors}
\end{figure}

From numerical tests we observed that both the third order schemes RK23 and Kutta's method behave like a second order scheme as well as the Midpoint scheme (as expected); this is not surprising since the drop evolution is coupled together with the second order IMEX scheme for the surfactant. Once we introduce adaptivity in the time-step selection based on both surfactant and drop evolution, we can observe that the \textit{conservation error} is a good choice for adaptivity only when using the Midpoint scheme, whilst for the other two schemes the error is not below the prescribed tolerance (fig. \ref{fig:tol_errors}). In the same figure we can see that when IMEX1 is used for the adaptivity in the surfactant, the error is lower, even too low for the drop evolution (fig. \ref{fig:tol_errors_drop}) when using RK23 or Kutta's method; in these two cases, the number of Stokes evaluations (which represents the main cost) is the same for the two schemes, and lower for the Midpoint scheme (fig. \ref{fig:tol_stokes_IMEX1}). The Midpoint scheme is the only one which respects the prescribed tolerance in the drop and in the surfactant using both IMEX1 and the conservation error. Moreover, since the adaptivity with the conservation error needs no extra evaluation, the Midpoint rule is the best choice. Note that for these simulations we used a rather large expansion order $p$ in order to make sure that the relative error is independent from the spatial resolution; the error flattens out when the tolerance decrease (Fig. \ref{fig:tol_errors}) due to the fixed imposed maximum step-size.

\begin{algorithm}[h!]
  \caption{Coupling the drop/surfactant evolution in time}
  \vspace{0.3cm}
  \begin{algorithmic}
\While{$t<T_{max}$}
\begin{itemize}
\item Given $\textbf{x}^{(t)},\Gamma^{(t)}$, compute velocity $\textbf{u}^{(t)}$ by solving eq. (\ref{eq:BIM_galerkin}).
\end{itemize}
\Large
$$\Downarrow$$
\normalsize
\begin{itemize}
\item $\textbf{x}^{(t+dt/2)}=\textbf{x}^{(t)}+\frac{dt}{2}\textbf{u}^{(t)};$
\item Compute $\Gamma^{(t+dt/2)}$ by solving the system in (\ref{eq:systemforgammatpdt2});
\item Given $\textbf{x}^{(t+dt/2)},\Gamma^{(t+dt/2)}$, compute velocity $\textbf{u}^{(t+dt/2)}$ by solving eq. (\ref{eq:BIM_galerkin}).
\end{itemize}
\Large
$$\Downarrow$$
\normalsize
\begin{itemize}
\item $\textbf{x}^{(t+dt)}=\textbf{x}^{(t+dt/2)}+dt\textbf{u}^{(t+dt/2)};$
\item Compute $\Gamma^{(t+dt)}$ using eq. (\ref{eq:gammatpiudt}).
\end{itemize}

\\

\begin{tikzpicture}

\draw [gray] (-2,-9.3) -- (10,-9.3);


\draw (-2.5,-14) rectangle (5,-10);
\node [right] at (0,-10.5) {\textit{Adaptivity in time}};
\node [right] at (-2.3,-11) {$err=max(err_{drop},err_{surfactant})$;};
\node [right] at (-2.3,-11.5) {\textbf{if} $err<tol$ \textbf{then} };
\node [right] at (-2,-12) {\textit{Reparametrization $\rightarrow$ Algorithm \ref{alg:pseudo_rep}};};
\node [right] at (-2,-12.5) {t=t+dt;};
\node [right] at (-2.3,-13) {\textbf{end if}};
\node [right] at (-2.3,-13.6) {$dt=dt(0.9 \frac{tol}{err})^{1/2}$};
\end{tikzpicture}
\\

 \EndWhile

   \end{algorithmic}
       \label{alg:pseudo_code}
  \end{algorithm}

 \vspace{-0.1cm}
 \vspace{-0.8cm}

  \clearpage
\newpage

\section{Reparametrization}
\label{reparametrization}
%
%
%
\begin{figure}[h!]
    \centering
 \begin{subfigure}[b]{0.3\textwidth}
         \includegraphics[width=60mm]{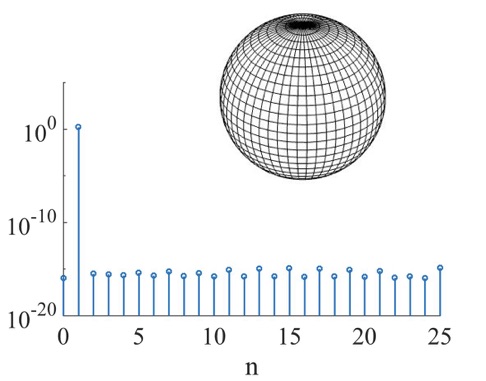}
         \caption{}
        \label{fig:sphere_nondist}
    \end{subfigure}
    \hspace{2cm}
      \begin{subfigure}[b]{0.3\textwidth}
         \includegraphics[width=60mm]{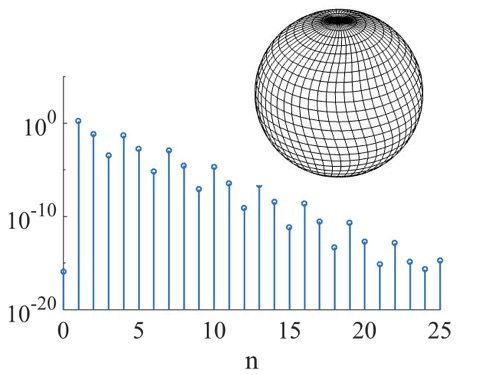}
        \caption{}
        \label{fig:sphere_dist}
    \end{subfigure}
 
    \caption{(a) Standard and (b) distorted grids on a unit sphere with relative spectra $E_n(\textbf{x}) := (\sum_{m=-n}^n |\textbf{x}_{n}^{m}|^2)^{1/2}$.}
  \label{fig:distorted}

\end{figure}

The spherical harmonics representation can be seen as a mapping from a Gaussian grid to the 3D surface. This mapping is not unique for a given surface (see Figure \ref{fig:distorted}) and the optimal map\footnote{Optimal in the sense of modes needed. For the easiest example of the sphere we know that the optimal choice is the one that needs only one mode, see Figure \ref{fig:sphere_nondist}.} is not known a priori for a general surface. When performing numerical simulations, significant distortions of the point distributions on the drop
surfaces may arise, especially for long simulations. This can
introduce unresolvable high-frequency components
(Fig. \ref{fig:sphere_dist}), aliasing errors and numerical
instability. For this reason we need to introduce a reparametrization
algorithm which seeks to minimize a spectral energy, optimizing the
grid point distribution with respect to the spherical harmonic
representation. 
We will begin by describing the approach used in
\cite{6,1}. Here, the authors introduce a smooth function
$F:\mathbb{R}^3 \rightarrow \mathbb{R}$, which is the implicit
representation of a general surface $\gamma$, and a quality measure
$E(\textbf{x}):\textit{X} \rightarrow \mathbb{R}$, where $\textit{X}$
is a space of sufficiently smooth functions on $\mathbb{S}^2$. It is
shown how the solution of the optimization problem 
\begin{equation}
\label{eq:opt_form_rep}
\min_{\textbf{x}\in \textit{X}} \{ E(\textbf{x}(\theta,\phi))\}, \ \text{subject to} \ F(\textbf{x}(\theta,\phi))=0,
\end{equation}
is given by
\begin{equation}
\label{eq:rep_prob_stat1}
\textbf{g}:=(I-\textbf{n}(\textbf{x}) \otimes \textbf{n}(\textbf{x}))
\nabla E(\textbf{x}) =0, \ F(\textbf{x})=0. 
\end{equation}
The authors \cite{6} suggest to use the quality measure (or energy) 
\begin{equation}
\label{eq:energy_rep}
E(\textbf{x})= \sum_{n=0}^p  \sum_{m=-n}^n a_{nm}|\textbf{x}_n^m|^2,
\end{equation} 
where $a_{nm}$ is some attenuation coefficient and $\textbf{x}_n^m$
are the coefficients of the spherical harmonic expansion of the
position vector $\textbf{x}$.

Eq. (\ref{eq:rep_prob_stat1})  can be solved by introducing a fictitious
time (or psuedo-time) $\tau$ and looking for stationary solution of 
\begin{equation}
\label{eq:pseudo_rep}
\begin{cases}
\frac{\partial \textbf{x}}{\partial \tau}+(I-\textbf{n}(\textbf{x})
\otimes \textbf{n}(\textbf{x})) \nabla E(\textbf{x}) =0, \\
\textbf{x}(0)=\textbf{x}_0.
\end{cases}
\end{equation}
A maximum number of iterations $i_{max}$ and a tolerance
$\epsilon$ are set for the breaking criteria, where $\epsilon$ is a
parameter to control the change in the distribution of points (see the original references for more details). Note that the condition $F(\textbf{x})=0$ is embedded in the continuous formulation. In \cite{6,1}, eq. (\ref{eq:pseudo_rep}) is discretized using the forward
Euler scheme and the discrete surface points are marched in pseudo
time. Here, the points are not constrained to stay on the surface, which introduces an error that depends on the pseudo time step size.

The aim of the reparametrization is to penalize the high frequency
components, i.e. the attenuation coefficients $a_{nm}$ should be small
for low frequencies and large for the highest ones. We can use a
\textsl{Perfect Low-Pass Filter}, where $a_{nm}=0$ if $n<n_{cutoff}$
and $a_{nm}=1$ otherwise, or we can also use a smoother filter that
increases with frequencies: $a_{nm}=0$ if $n<n_{cutoff}$ and
$a_{nm}=\frac{n}{p}$ otherwise. 
In \cite{6,1} such reparametrization with a \textsl{Perfect Low-Pass
  Filter} and a fixed $n_{cutoff}$ ($n_{cutoff}=\frac{p}{k_{cut}}$,
$k_{cut}=2,3$) is used for the simulation of vesicles where it gives
good results. 

The main difference between vesicles and drops in this context, is
that for the first case the local surface inextensibility avoids the
same degree of distortions of the point distribution that can be
seen for drops.
We encounter a number of difficulties: 
\begin{enumerate}[label=(\alph*)]
\item \label{rep1} We note that the above mentioned
  reparametrization will not ensure volume conservation; 
\item \label{rep2}  Using a fixed $n_{cutoff}$ is not working well in our simulations, neither the Perfect Low-Pass Filter nor the smoother one; 
\item \label{rep3} The current procedure offers no spectrally accurate
  way to redefine the surfactant concentration on the reparameterized
  surface. 
\end{enumerate}
As for the first remark, it is certainly true also for vesicles that
the volume conservation error will depend on the pseudo time step
size, but it should be more noticeable for drops where distortions are
larger and more reinitialization is needed. 

In the following, we will introduce a new reparameterization method that resolves issues \ref{rep1} and \ref{rep3}. 
We will however start by introducing an adaptive $n_{cutoff}$ to overcome the
problem \ref{rep2}, and show how this
change by itself improves also the old reparameterization procedure. 

\begin{figure}[htbp]
	\centering
			\includegraphics[width=90mm]{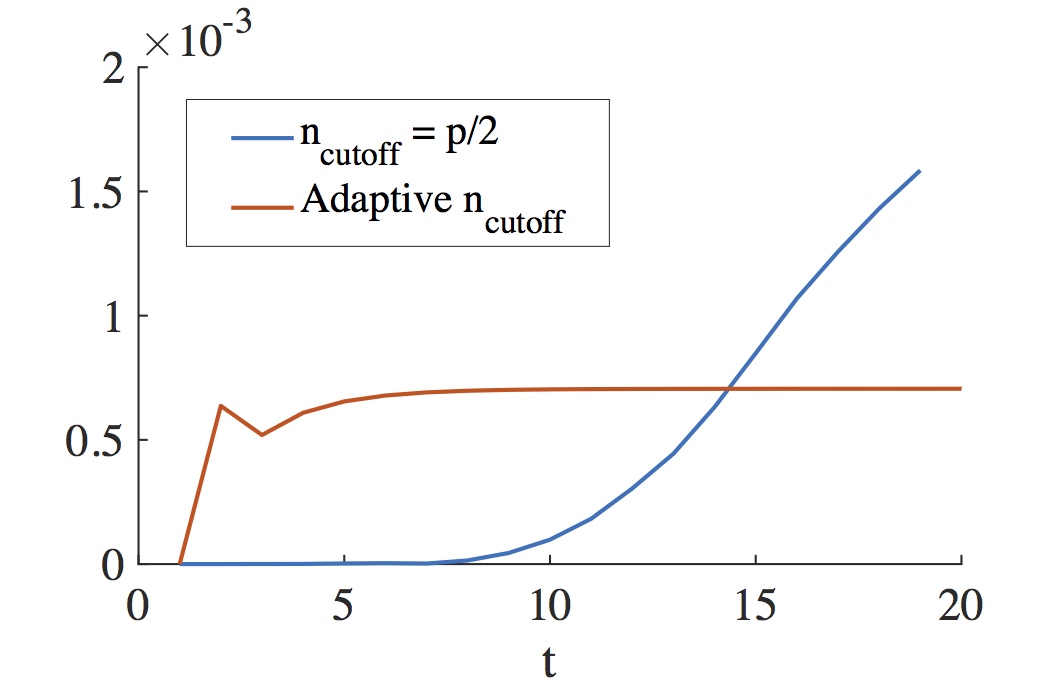}
	\caption{Comparing the quality measure $E$ defined in (\ref{eq:energy_rep}) using adaptive and non-adaptive
          $n_{cutoff}$. This quantity is not stable in time when using the fixed $n_{cutoff}=p/2$. Instead, using the adaptive $n_{cutoff}$ the steady state is reached in good agreement with experimental results \cite{12}.}
	\label{fig:energy_rep}
\end{figure}

The adaptive $n_{cutoff}$ will depend on the geometry of the surface;
in particular we fix the percentage $P_{cutoff}$ of the modes that we
want to penalize and then we compute $n_{cutoff}$ as:
\begin{equation}
\label{eq:ncutoff_adaptive}
n_{cutoff}=min_l \Big \{ l \in \mathbb{N}, 1<l \leq p : \frac{N_2(l)}{N_1}<P_{cutoff} \Big\}
\end{equation}
where $N_2(l)=\sum_{n=l}^p\sum_{m=-n}^n |\textbf{x}_n^m|$ and
$N_1=N_2(1)$. This adaptive choice can much better determine what
frequencies need to be cut for a given geometry. A basic example is a
distorted grid on a sphere (as in Fig. \ref{fig:sphere_dist}): in this
case we know that only one mode is needed for the optimal distribution
of points (Fig. \ref{fig:sphere_nondist}), and we reach this state
using the adaptive $n_{cutoff}$ but not using a fixed $n_{cutoff}$ (in
which case only modes higher than $p/2$ would be filtered out).
Numerical experiments show that the new procedure works better also for
more complicated problems, like drops immersed in a linear flow that
reach a steady state. In this case the drop shape is steady but due to
tangential velocities on the drop surface the point distribution is
continuously changing and the reparametrization plays a fundamental
role. To show this, we consider a clean drop immersed in a flow
generated in the four-roll mill
\[\nabla u=\frac{1}{2}G
  \left[ {\begin{array}{ccc}
   1+\alpha & 1-\alpha & 0 \\  -1+\alpha & -1-\alpha & 0 \\ 0 & 0 & 0 \ \end{array} } \right],
\]
where $\alpha$ is a parameter which specifies the relative strength of
the strain rate and vorticity in the flow. The parameters which affect
the evolution of the drop are set to be $Ca=0.0718$, $\lambda=0.118$
and $\alpha=0.6$. Under these conditions, the drop reaches a steady
state, meaning that the energy defined in eq. (\ref{eq:energy_rep})
should reach a constant value. In Fig. \ref{fig:energy_rep} this
expected behavior can be seen using the adaptive $n_{cutoff}$, but for
a fixed $n_{cutoff}=p/2$ the energy continues to increase and yields
an unstable simulation.  Defining the deformation number as:
\begin{equation}
\label{eq:def_num}
D=\frac{L-B}{L+B}
\end{equation} 
where $L$ and $B$ are the drop length and breadth in the flow plane
respectively, for the simulation with the adaptive $n_{cutoff}$ we
find a value of $D=0.117$, in good agreement with the experimental
results presented in \cite{12}.
Since this adaptive filter works well it will be used also with the
new reparametrization procedure.

As mentioned in \ref{rep3}, we have a quantity (the surfactant
concentration) that is defined on the surface.  Before starting the
reinitialization procedure, we have the initial expansions
$\mathbf{x}_0(\theta,\phi)$ for the drop surface and
$\Gamma_0(\theta,\phi)$ for the surfactant concentration.  The discrete
points on the surface are given by $\mathbf{x}(\theta,\phi)$ evaluated in the
$\theta$ and $\phi$ values as given in (\ref{eq:surface_ev}).  In the old
reinitialization procedure, we move the discrete points, and for this reason we will refer to this procedure as to the \textit{point-reparametrization}.  After they
have been moved, we can compute new coefficients for this discrete set
of points using the discrete transform (\ref{eq:back_transf_discr}). We do not however know what
$\theta$ and $\phi$ values the new point locations correspond to in
the initial expansion.  If we knew that, we could compute the
corresponding $\Gamma$ values by simply evaluating $\Gamma_0(\theta,\phi)$
at these $(\theta,\phi)$ values. Hence, it is the aim for the new
reinitialization procedure to find these values of the angles, and for this reason we will refer to the new procedure as to the \textit{angle-reparametrization}. As an additional
benefit, knowing these $\theta$ and $\phi$ values, we can evaluate the
discrete points on the surface from the initial expansion
$\mathbf{x}_0(\theta,\phi)$, which yields a spectrally accurate volume
conservation.

\begin{figure}[htbp]
	\centering
			\includegraphics[width=150mm]{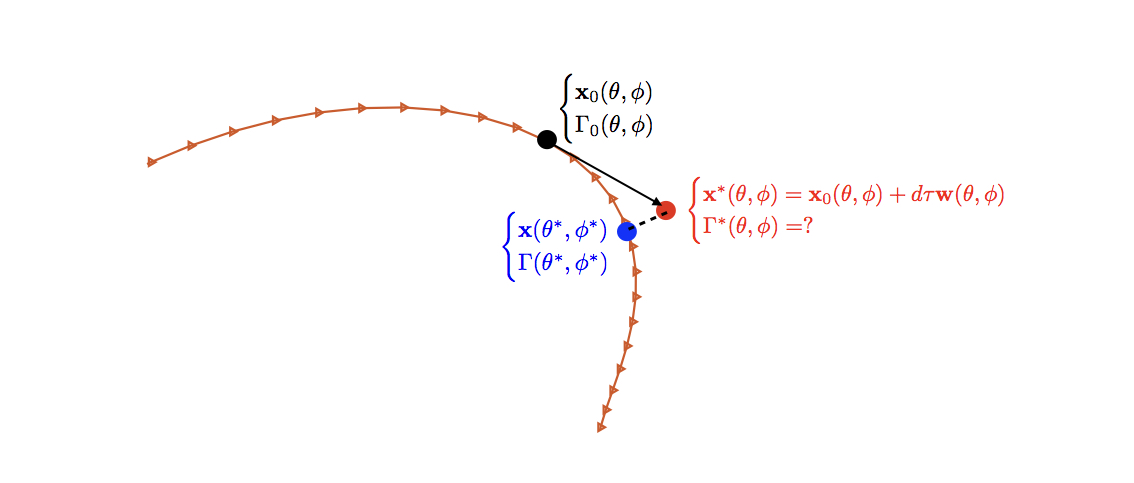}
	\caption{Sketch of the problem for the reparametrization.}
	\label{fig:sketch_rep}
	
\end{figure}

Fig. \ref{fig:sketch_rep} shows a 2D-sketch of the problem, also
indicating that the discrete points are not constrained to stay on the
surface in the \textit{point-reparametrization} procedure. For the new procedure,  
we instead need to find the new angle vector $(\theta^*,\phi^*)$ in
order to evaluate both the point location and the surfactant
concentration on the surface. To do this we compute a parameter pseudo-velocity $(\frac{\partial \theta}{\partial \tau},\frac{\partial \phi}{\partial \tau})$ for evolving $(\theta,\phi)$ in pseudo-time. From eq. (\ref{eq:pseudo_rep}) we have $\frac{\partial \textbf{x}}{\partial \tau}=\textbf{w}$, where the pseudo-velocity $\textbf{w}$ is given by
$$\textbf{w}(\theta,\phi)=[u_1(\theta,\phi),u_2(\theta,\phi),u_3(\theta,\phi)]=-(I-\textbf{n}(\textbf{x})
\otimes \textbf{n}(\textbf{x})) \nabla E(\textbf{x}).$$
By the chain rule, we have 
\begin{equation}
\frac{d\textbf{x}}{d \tau}=\frac{\partial \textbf{x}}{\partial \theta} \frac{\partial \theta}{\partial \tau}+\frac{\partial \textbf{x}}{\partial \phi} \frac{\partial \phi}{\partial \tau}.
\end{equation}
In order to find the parameter pseudo-velocity $(\frac{\partial \theta}{\partial \tau},\frac{\partial \phi}{\partial \tau})$, we consider the two projections, where $\langle\cdot,\cdot\rangle$ denotes the scalar product:
\begin{equation*}
\begin{cases}
            \langle\frac{d\textbf{x}}{d \tau},\frac{\partial \textbf{x}}{\partial \theta}\rangle=\langle\frac{\partial \textbf{x}}{\partial \theta} \frac{\partial \theta}{\partial \tau}+\frac{\partial \textbf{x}}{\partial \phi} \frac{\partial \phi}{\partial \tau},\frac{\partial \textbf{x}}{\partial \theta}\rangle\\
\langle\frac{d\textbf{x}}{d \tau},\frac{\partial \textbf{x}}{\partial \phi}\rangle=\langle\frac{\partial \textbf{x}}{\partial \theta} \frac{\partial \theta}{\partial \tau}+\frac{\partial \textbf{x}}{\partial \phi} \frac{\partial \phi}{\partial \tau},\frac{\partial \textbf{x}}{\partial \phi}\rangle

       \end{cases} \rightarrow
			\quad
 \begin{cases}
            \langle\textbf{w},\frac{\partial \textbf{x}}{\partial \theta}\rangle=\frac{\partial \theta}{\partial \tau}|\frac{\partial \textbf{x}}{\partial \theta}|^2 +\frac{\partial \phi}{\partial \tau}\langle\frac{\partial \textbf{x}}{\partial \phi} ,\frac{\partial \textbf{x}}{\partial \theta}\rangle\\
\langle\textbf{w},\frac{\partial \textbf{x}}{\partial \phi}\rangle=\frac{\partial \theta}{\partial \tau}\langle\frac{\partial \textbf{x}}{\partial \phi} ,\frac{\partial \textbf{x}}{\partial \theta}\rangle +\frac{\partial \phi}{\partial \tau}|\frac{\partial \textbf{x}}{\partial \phi}|^2
       \end{cases}
\end{equation*}
which can be stated in compact form as
\begin{equation}
\label{eq:system_rep}
A\alpha=b,
\end{equation}
where
\begin{equation*}
A=\left[\begin{matrix} |\frac{\partial \textbf{x}}{\partial \theta}|^2 & \langle\frac{\partial \textbf{x}}{\partial \phi} ,\frac{\partial \textbf{x}}{\partial \theta}\rangle \\ \langle\frac{\partial \textbf{x}}{\partial \phi} ,\frac{\partial \textbf{x}}{\partial \theta}\rangle  & |\frac{\partial \textbf{x}}{\partial \phi}|^2
\end{matrix}\right], \ 
\alpha=\left(
\begin{array}{c}
\frac{\partial \theta}{\partial \tau}\\
\frac{\partial \phi}{\partial \tau}\\
\end{array}
\right), \ b= \left(
\begin{array}{c}
\langle\textbf{w},\frac{\partial \textbf{x}}{\partial \theta}\rangle\\
\langle\textbf{w},\frac{\partial \textbf{x}}{\partial \phi}\rangle\\
\end{array}
\right).
\end{equation*}
We can explicitly write
\begin{equation}
\label{eq:system_rep_sol}
\alpha=A^{-1}b,
\end{equation}
where 
$$A^{-1}=\frac{1}{|\frac{\partial \textbf{x}}{\partial \theta}|^2|\frac{\partial \textbf{x}}{\partial \phi}|^2 -\langle\frac{\partial \textbf{x}}{\partial \phi} ,\frac{\partial \textbf{x}}{\partial \theta}\rangle^2}\left[\begin{matrix} |\frac{\partial \textbf{x}}{\partial \phi}|^2 & -\langle\frac{\partial \textbf{x}}{\partial \phi} ,\frac{\partial \textbf{x}}{\partial \theta}\rangle \\-\langle\frac{\partial \textbf{x}}{\partial \phi} ,\frac{\partial \textbf{x}}{\partial \theta} \rangle & |\frac{\partial \textbf{x}}{\partial \theta}|^2
\end{matrix}\right].$$

Once we solve the system for the derivatives $(\frac{\partial \theta}{\partial \tau}, \frac{\partial \phi}{\partial \tau})$, we can compute the new parameters $(\theta^*,\phi^*)$ using a forward Euler step. At this point the new position is given by eq. (\ref{eq:surface_ev}). One extra forward and one extra backward transforms are needed in order to
keep the original grid, equidistant in $\phi$ and Gauss-Legendre points in $\theta$:
\begin{equation} 
\label{eq:rep_trasf1} 
\textbf{x}(\theta^*,\phi^*)=\sum_{n=0}^p \sum_{m=-n}^n
\textbf{x}_{(0),n}^m
Y_n^m(\theta^*,\phi^*)=:\tilde{\textbf{x}}(\theta,\phi)=\sum_{n=0}^p \sum_{m=-n}^n
\tilde{\textbf{x}}_{n}^m
Y_n^m(\theta,\phi)
\end{equation} 
\begin{equation}
 \label{eq:rep_trasf2}
\Gamma(\theta^*,\phi^*)=\sum_{n=0}^p \sum_{m=-n}^n \Gamma_{(0),n}^m Y_n^m(\theta^*,\phi^*)=:\tilde{\Gamma}(\theta,\phi)=\sum_{n=0}^p \sum_{m=-n}^n
\tilde{\Gamma}_{n}^m
Y_n^m(\theta,\phi),
\end{equation}
where $\tilde{\textbf{x}}_{n}^m$ and $\tilde{\Gamma}_{n}^m$ are given by eq. (\ref{eq:back_transf}).

\begin{figure}[htbp]
	\centering
	\begin{subfigure}[b]{0.4\textwidth}
		\includegraphics[width=65mm]{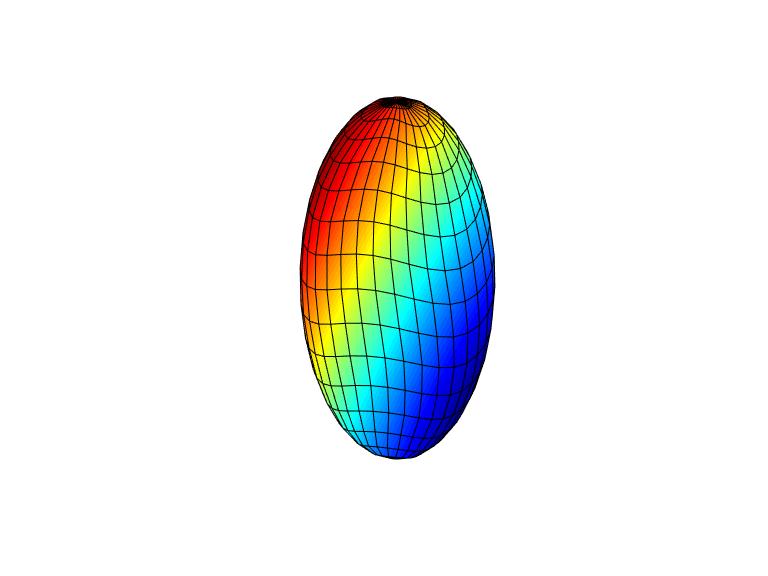}
		\caption{}
		\label{fig:dist_ell_tp}
		\end{subfigure}
		\quad
		\begin{subfigure}[b]{0.4\textwidth}
		\includegraphics[width=65mm]
{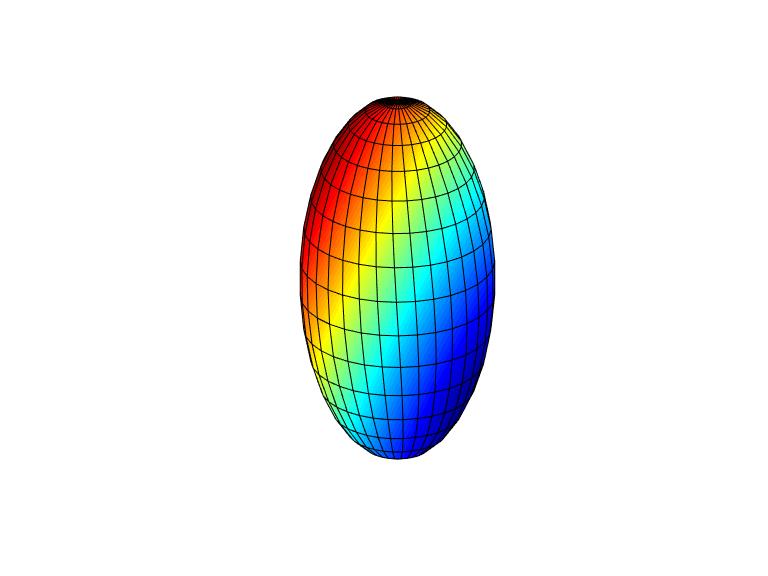}
\caption{}
\label{fig:dist_ell_tp_after}
		\end{subfigure}
\caption{The distorted ellipsoid with $\Gamma=2-\frac{1}{2}x_1+x_2+\frac{1}{2}x_3$, (a) pre and (b) after reparametrization. Colors denote the surfactant concentration.}
  \label{fig:test_surf_rep_ell}
\end{figure}

 This whole procedure is done for all discrete points on the surface in each pseudo time step, and it is summarized in the pseudo-code Algorithm \ref{alg:pseudo_rep}.\\
\begin{algorithm}
  \caption{Reparametrization}
  \vspace{0.3cm}
  \begin{flushleft}
   Given an initial configuration for the drop surface and for the surfactant concentration we upsample by a factor $U_{rep}$;  on this finer grid the initial drop surface and initial surfactant concentration are given by $[\textbf{x}_0(\theta,\phi),\Gamma_0(\theta,\phi)]$ with spectra $[\textbf{x}_{(0),n}^m,\Gamma_{(0),n}^m]$. Set $\tilde{\textbf{x}}(\theta,\phi)=\textbf{x}_0(\theta,\phi)$ and $\tilde{\textbf{x}}_n^m=\textbf{x}_{(0),n}^m$
  \end{flushleft}
  \begin{algorithmic}
\While{++$i<i_{max}$}\\
\begin{itemize}
\item Compute the energy using eq. (\ref{eq:energy_rep}) and adaptive $n_{cutoff}$ computed through eq. (\ref{eq:ncutoff_adaptive});
\item Use the computed energy for obtaining the pseudo-velocity $\frac{\partial \textbf{x}}{\partial \tau}=\textbf{w}(\theta,\phi)$ through eq. (\ref{eq:pseudo_rep});
\item Solve system (\ref{eq:system_rep}) to find the parameter pseudo-velocity $\alpha=(\frac{\partial \theta}{\partial \tau}, \frac{\partial \phi}{\partial \tau})$;
\item Compute the new parameters $(\theta^*,\phi^*)$ using a forward Euler step;
\item Update the position vector using the initial expansion coefficients, eq. (\ref{eq:rep_trasf1})
\item Compute the updated expansion coefficients $\tilde{\textbf{x}}_n^m$ via (\ref{eq:back_transf_discr})
\item Compute the normal vector and the gradient of the energy (\ref{eq:energy_rep}) corresponding to the new position in order to update $\textbf{g}$ via eq. (\ref{eq:rep_prob_stat1})
\end{itemize}
\If {$\|\textbf{g}\|<\epsilon$} 
break
\EndIf

 \EndWhile
 \begin{itemize}
 \item Compute the new surfactant concentration using the initial expansion coefficients, eq. (\ref{eq:rep_trasf2});
\item Compute the updated expansion coefficients $\tilde{\Gamma}_n^m$ via (\ref{eq:back_transf_discr}).
 \end{itemize}
   \end{algorithmic}
       \label{alg:pseudo_rep}
  \end{algorithm}
To test the new reparametrization, we intentionally
distort\footnote{We apply an artificial tangential velocity given by
  $\textbf{v}=3kcos(3\theta)\textbf{t}_{\phi}+2kcos(3\phi)\textbf{t}_{\theta}$
  until the computational time $T=3$, with $k=0.01$} 
the grid of an ellipsoid with aspect ratio (1:2) as shown in Fig.
\ref{fig:dist_ell_tp} and let the new procedure run. 

\begin{figure}[htbp]
	\centering
	\begin{subfigure}[b]{0.42\textwidth}
		\includegraphics[width=70mm]{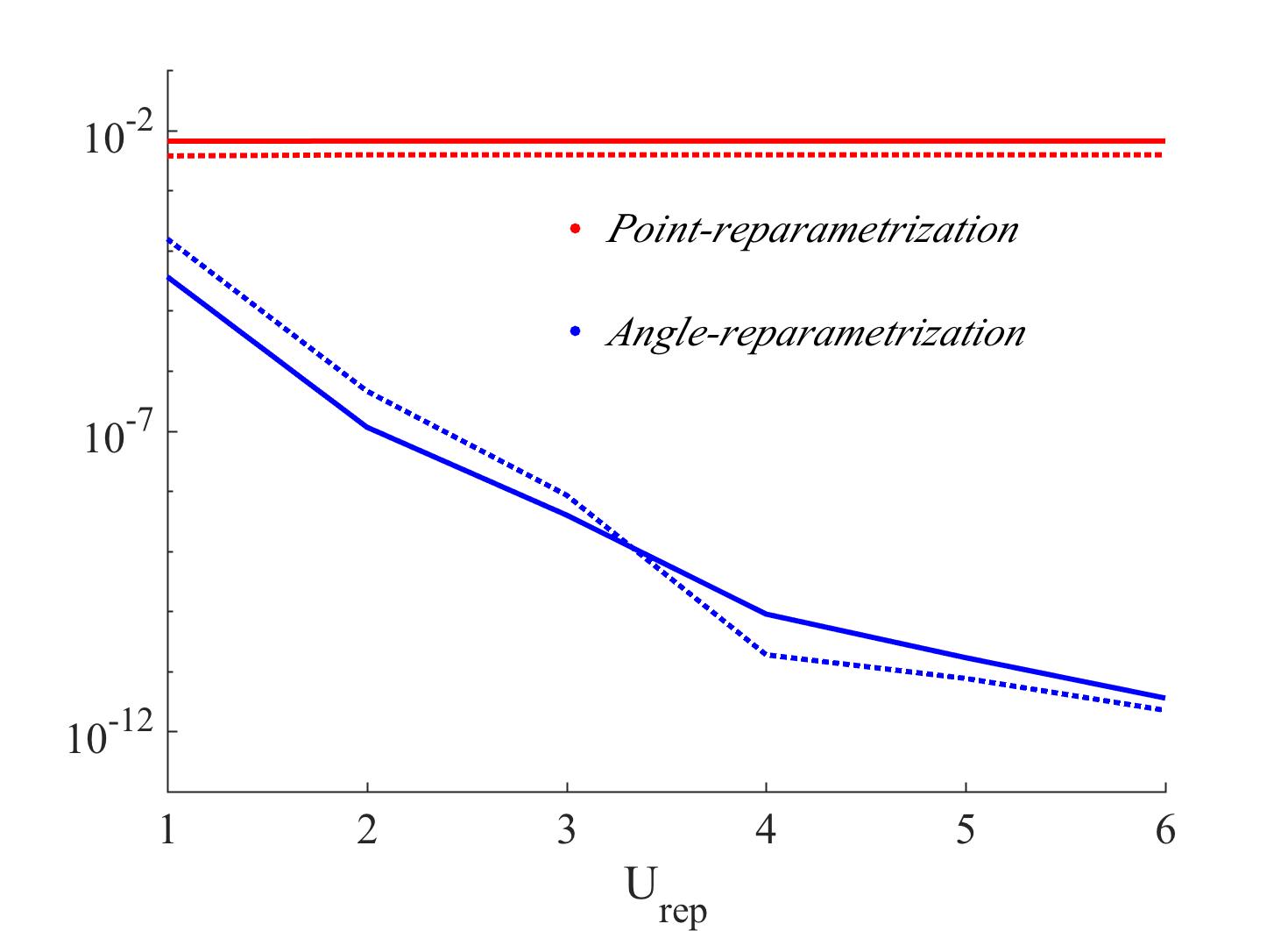}
		\caption{Area (dotted line) and volume (continuous line) relative errors vs the upsampling factor U$_{rep}$ with $d\tau=1$.}
		\label{fig:UvsError}
		\end{subfigure}
		\quad
		\begin{subfigure}[b]{0.42\textwidth}
		\includegraphics[width=70mm]
{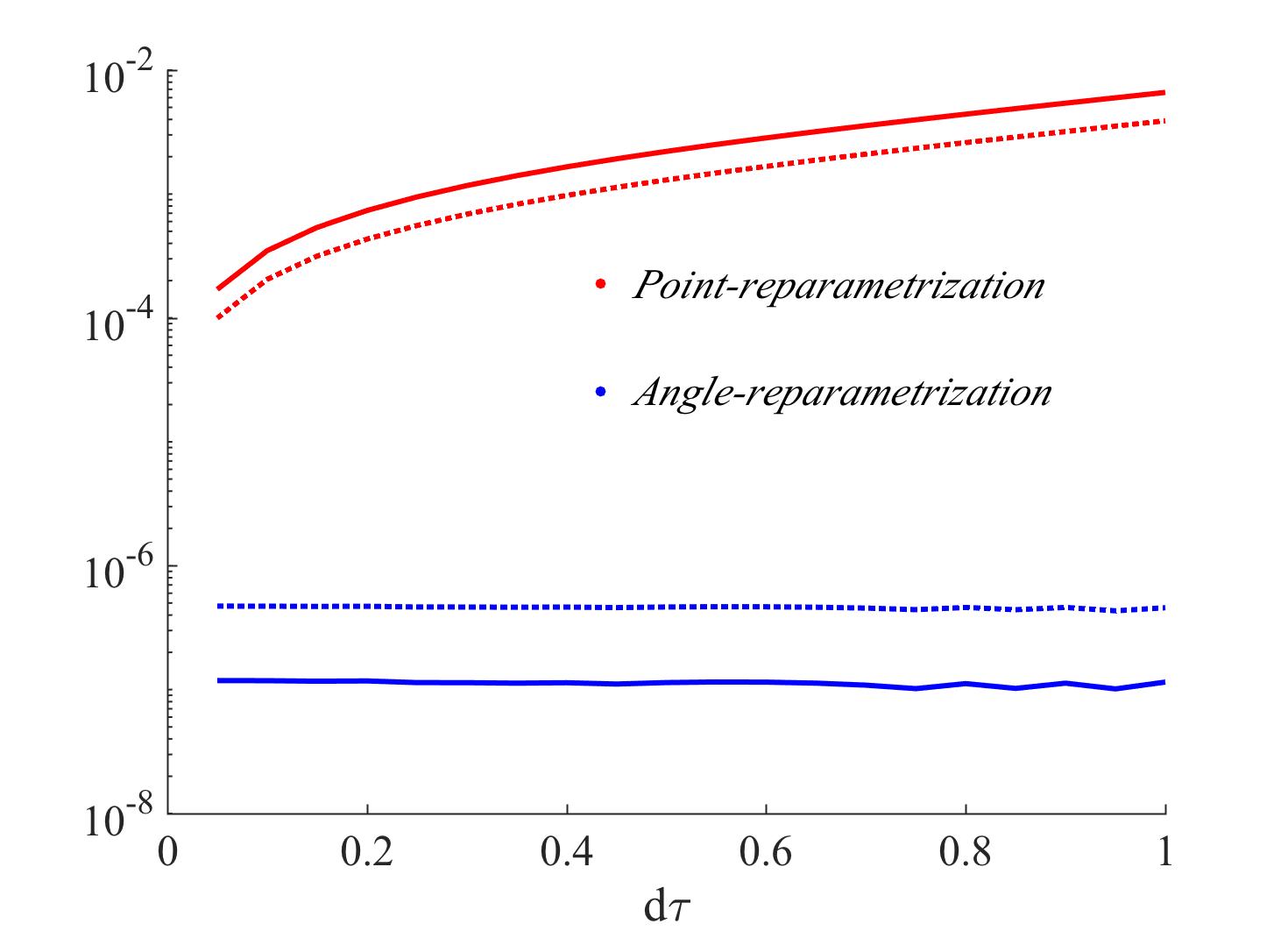}
\caption{Area (dotted line) and volume (continuous line) relative errors vs the pseudo time-step $d\tau$ with U$_{rep}=2$.}
\label{fig:dTvsError}
		\end{subfigure}
\caption{Comparing the \textit{point-reparametrization} and the new \textit{angle-reparametrization}, $p=7$, for the test case shown in Fig. \ref{fig:test_surf_rep_ell}.}
  \label{fig:test_repnew_ell_both_dt1_resp_P}
\end{figure}

The reparametrization is following only a fictitious velocity, so the
surfactant distribution should not change but appear the same as in
Fig.\ref{fig:dist_ell_tp} also after the reparametrization (Fig.
\ref{fig:dist_ell_tp_after}). In Fig.
\ref{fig:test_repnew_ell_both_dt1_resp_P} we compare the error in the
volume and the area conservations for the two methods of
reparametrization, using in each case the adaptive $n_{cutoff}$ and a
\textit{Perfect Low-Pass Filter}. We can see the spectral decay of the
error for the new reparametrization and, moreover, the
dependence/independence from the pseudo-time step $\tau$ can clearly
be seen (Fig. \ref{fig:test_repnew_ell_both_dt1_resp_P}): the error
for the \textit{point-reparametrization} is growing with this parameter, whilst
for the \textit{angle-reparametrization} it stays fixed and much smaller. This is because, with
the \textit{point-reparametrization}, we want to move the points tangential to
the surface, but we practically have to deal with a discrete system
and then a large pseudo-time step might bring the new point
distribution \textsl{far} from the surface, as shown in the sketch of
Fig. \ref{fig:sketch_rep}. With the new method, instead, we compute
the new point distribution and the new surfactant concentration by
evaluating a spherical transform, using the spherical harmonic
coefficients $\textbf{x}_{(0),n}^m$ and $\Gamma_{(0),n}^m$ of the
original surface/surfactant concentration, see eqs. (\ref{eq:rep_trasf1}-\ref{eq:rep_trasf2}). We are basically projecting the new position vector on the surface and computing the corresponding angle parameters; for this reason the size of the pseudo-time step is not a restriction anymore. Note also that the level of the error for the new reparametrization procedure is set by the upsampling rate of the reparametrization $U_{rep}$. This is shown in the example of Fig. \ref{fig:test_repnew_ell_both_dt1_resp_P}: if we compare the level of the error in Fig. \ref{fig:dTvsError} we can see that it corresponds to the error in Fig. \ref{fig:UvsError} for $U_{rep}=2$.

Finally we can check the surfactant concentration accuracy after the reparametrization using an analytical solution (Fig. \ref{fig:dist_ell_tp_after}):
$$\Gamma=2-\frac{1}{2}x_1+x_2+\frac{1}{2}x_3$$
and we found that the error is very low ($<10^{-15}$). This is obvious since we are just evaluating a spherical transform for a well-resolved geometry. We checked the reparametrization procedure for different geometries and surfactant distributions and results behave as for the previous case.

\section{Numerical experiments and validation}
\label{experiments}
In this section we will check our numerical method and perform some numerical experiments. A clean drop case has already been validated in Section \ref{reparametrization} comparing our results with the experiments by Bentley and Leal \cite{12}, so we proceeded by validating our code with standard tests for the surfactant convection-diffusion equation. In order to make the reader more familiar with the role of each single term acting in eq. (\ref{eq:surfactant0}), we describe the tests below:
\begin{enumerate}
\item \textit{Diffusion term}: in this case we do not consider the convective and the stretching term and, starting by a non-constant initial $\Gamma$, we can see the effect of the diffusion spreading the surfactant on the whole surface until it reaches a constant concentration, slower or faster depending on the P\'eclet number.
\item \textit{Stretching term}: we consider the uniform expansion of a spherical drop, radius $R(t)$, uniformly coated with insoluble surfactant concentration $\Gamma(t)$. The P\'eclet number is set to infinity so that only the stretching term is active (there is only the normal velocity component, so the convection term is zero). In this case we have the analytical solution $$\Gamma R^2(t)=\Gamma(0),$$ so that $\frac{\partial \Gamma}{\partial t}+2\frac{\Gamma}{R}\frac{\partial R}{\partial t}$=0, which is in agreement with equation (\ref{eq:surfactant0}), since the mean curvature of the sphere is $-\frac{1}{R}$ and then the stretching term is exactly $(\nabla_S \cdot n)(u \cdot n)=\frac{2}{R}\frac{\partial R}{\partial t}$.
\item \textit{Convective term:} for this case we consider a rotating sphere again with infinite P\'eclet number. In this case only the tangential velocity is non-zero and the solution is kept constant.
\end{enumerate}

At this point we are able to test the whole system, drop/surfactant for a single or multiple interacting drops. \\
\begin{figure}[htbp]
	\centering
	\begin{subfigure}[b]{0.42\textwidth}
	\includegraphics[width=70mm]{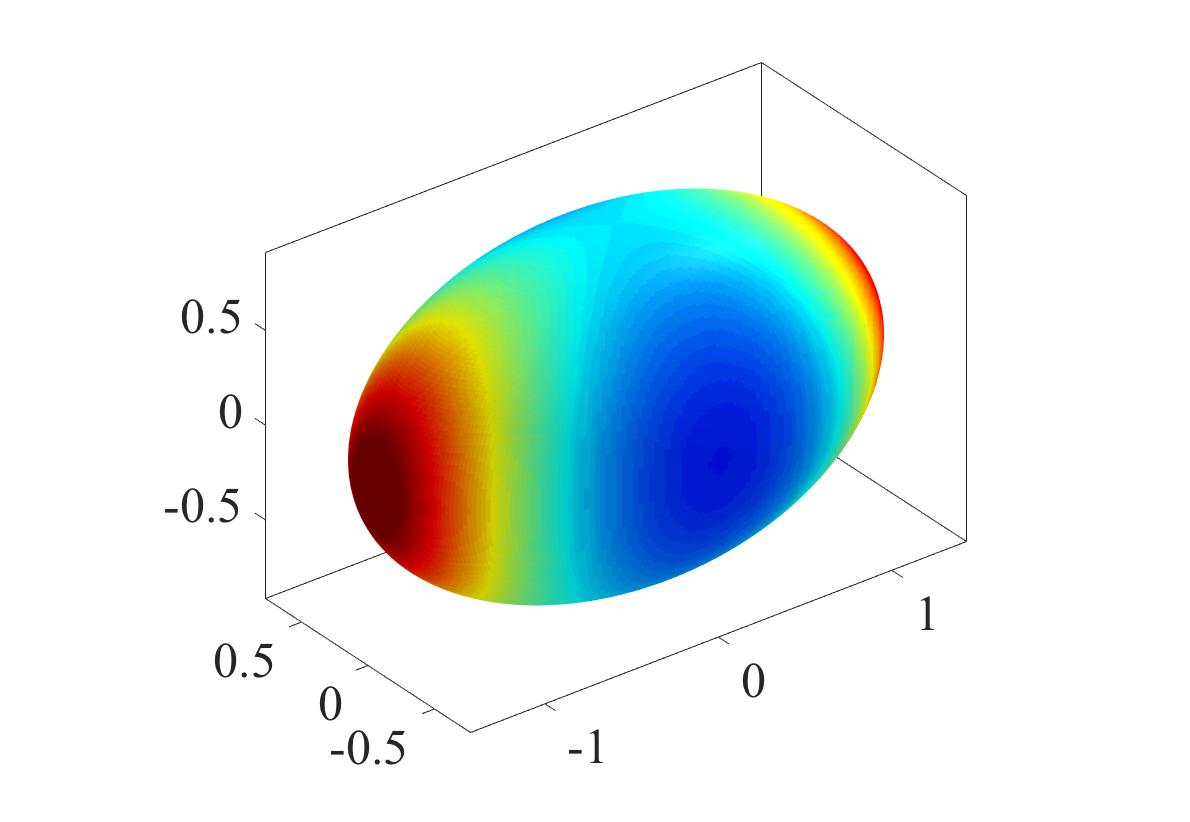}
	\caption{}
	\label{fig:test_ext3Dx036}
		\end{subfigure}
	\quad
		\begin{subfigure}[b]{0.42\textwidth}
		\includegraphics[width=70mm]{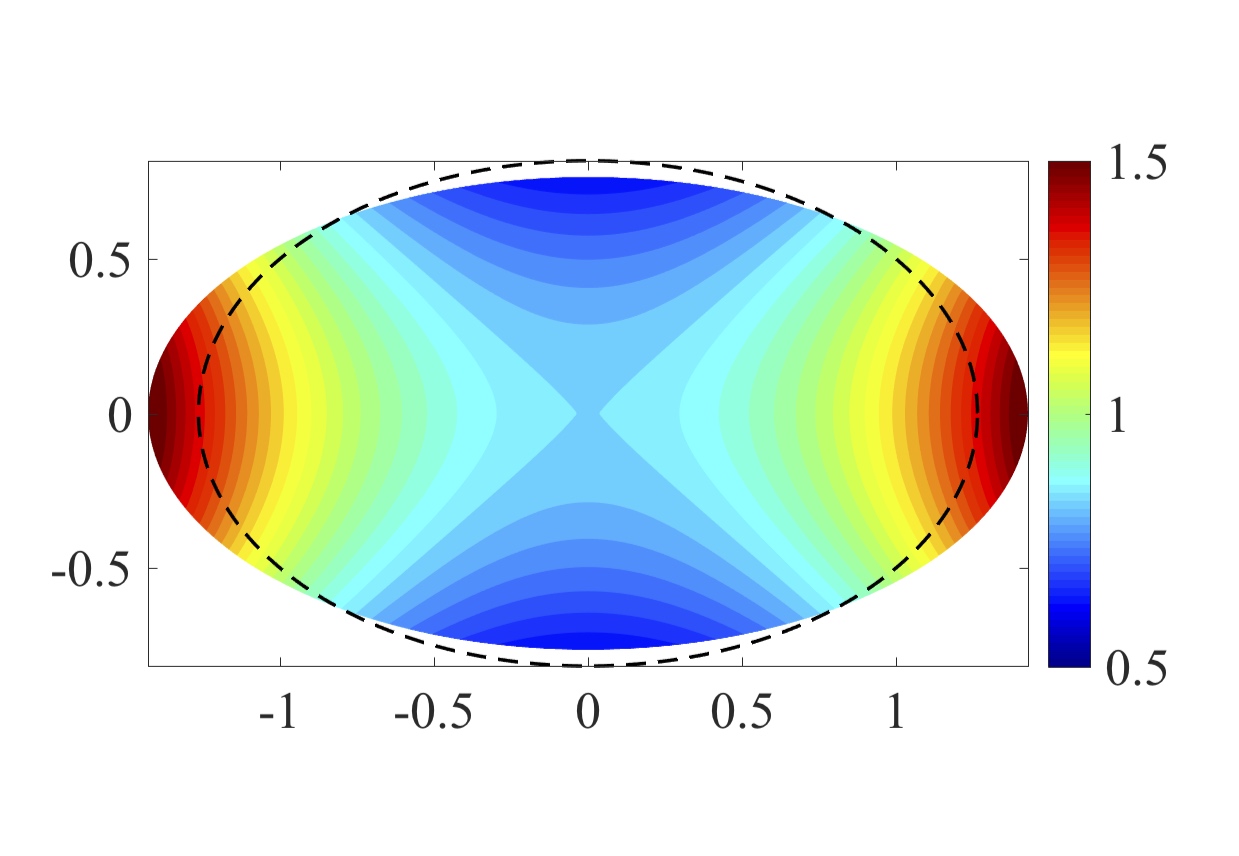}
		\caption{}
		\label{fig:test_ext_both}
		\end{subfigure}
	\caption{A single drop in an extensional flow with $E=0.35$, $Pe=11.8$, $Ca=0.1$, $\lambda=0.093$, $p=17$. (a) 3D surfactant-covered drop ($x_s=0.36$) at the steady state; (b) 2D projection to compare clean (black dotted line) and surfactant-covered (colorful) drops at the steady state.}
	\label{fig:test_hulips}
	
\end{figure}	
We start considering a drop with initial spherical shape subject to a planar extensional flow $\textbf{u}=(x,-y,0)$ immersed in a fluid with a different viscosity, $\lambda=0.093$. The specific value of $\lambda$ has been set by the experimental results \cite{22} and makes all the terms involved in our equation to be active: when $\lambda \neq 1$ we know that both the single and the double-layer play a fundamental role; moreover, when $\lambda$ is small and the surfactant coverage is big enough, the Marangoni term becomes significant \cite{18} and so the effect of the surfactant. We compare the behavior of a clean drop, with a surfactant-covered drop with surface-coverage (eq. \ref{eq:langmuir}) $x_s=0.36$, see Fig. \ref{fig:test_hulips}. We can see that the surfactant acts to lower the surface tension and consequently the drop deforms more. We found for these experiments the deformation numbers  (\ref{eq:def_num}) $D_{x_s=0} \approx 0.21$ and $D_{x_s=0.36}\approx 0.30$. These results are in good agreement with both the experimental results by Hu and Lips \cite{22} and with the numerical results by Bazhlekov et al.  \cite{18} (see Fig. 8 in \cite{18} where both experimental and numerical results are shown). We also performed a convergence test for this example: computing a reference solution ($p=29$) we look at the error for the drop position and surfactant concentration at a steady state, see Fig. \ref{fig:conv_plots}. This is a particularly interesting and complicated example to test: it's a time dependent problem, we are looking at a long time horizon and there is a strong tangential flow that distorts the point distribution on the drop even when the shape is steady. For this reason, the reparameterization is constantly activated and, since this procedure is not fixed when changing the expansion order $p$, it could ruin smooth error decay with increasing resolution; this is another reason why looking at the converge plot is so challenging. Nevertheless, we can see the spectral convergence of the whole method when the resolution is good enough.

\begin{figure}[htbp]
	\centering
		\hspace{-1cm}
		\includegraphics[width=80mm]{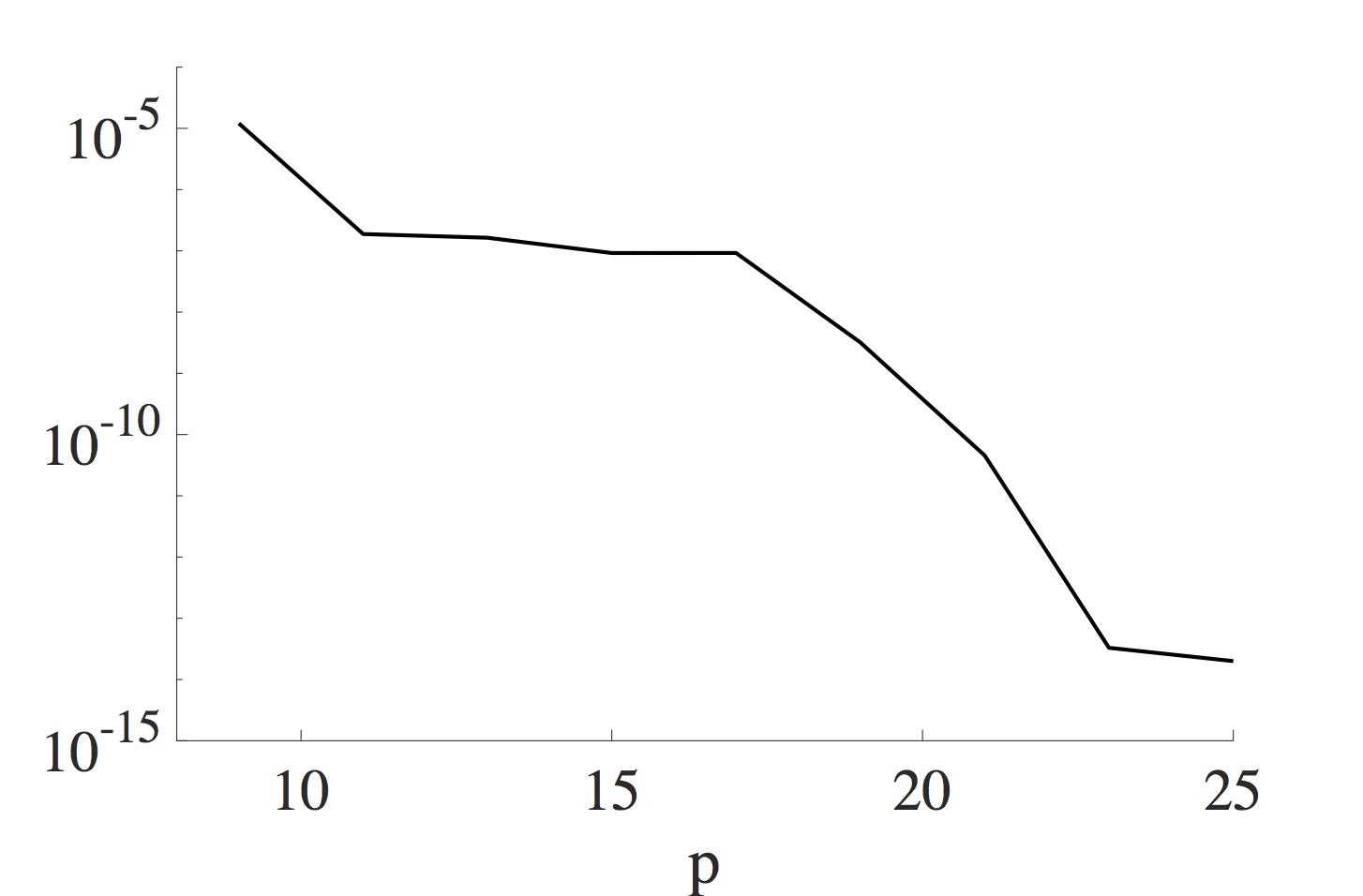}
		\quad
		\hspace{-1cm}
	\includegraphics[width=80mm]{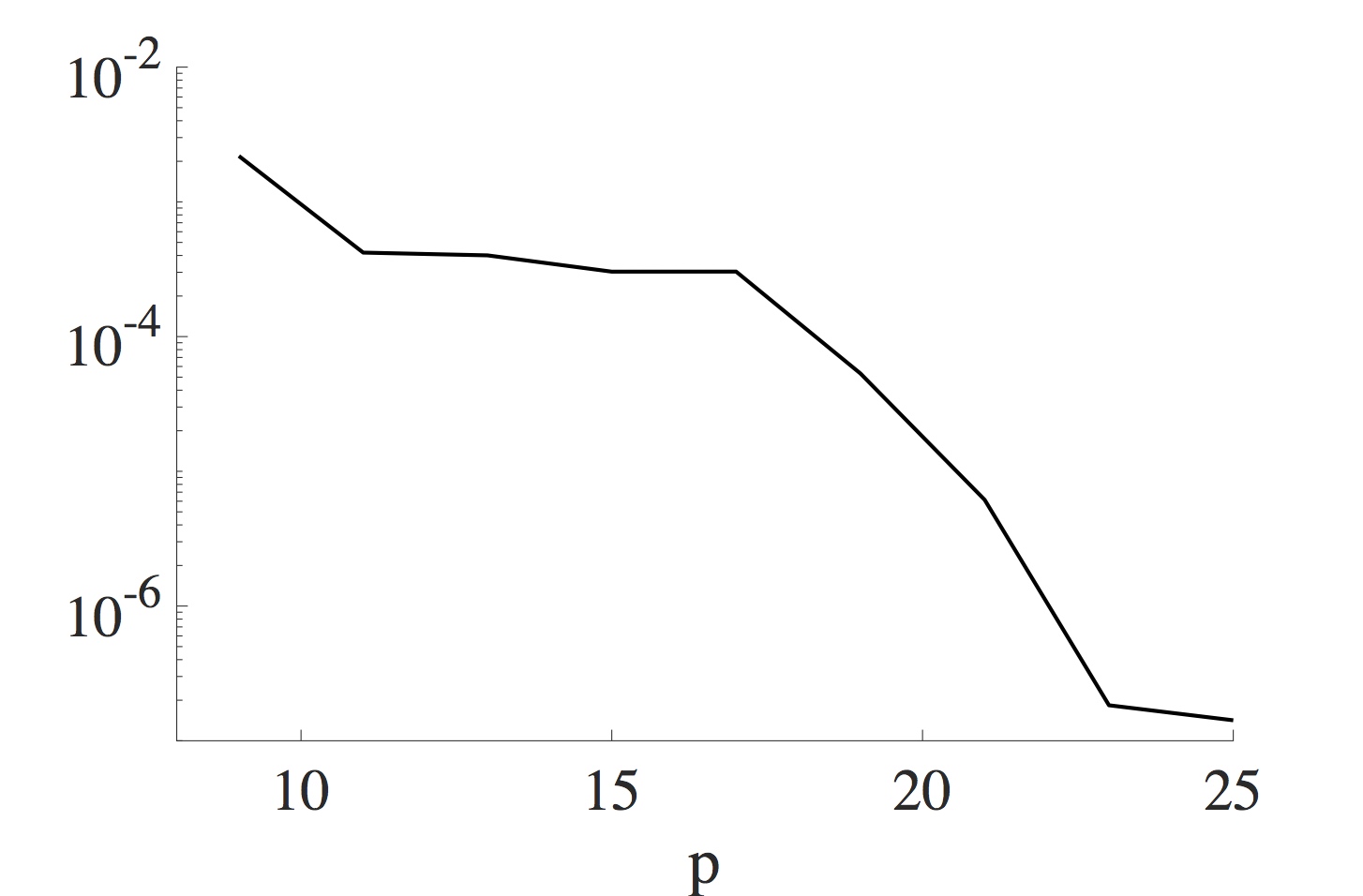}
	\caption{L2 norm of the error for the validation test of Fig. \ref{fig:test_hulips}. The error is computed for the spherical harmonics coefficients vector of the drop position (left) and of the surfactant concentration (right) using a reference solution obtained with $p=29$.}
	\label{fig:conv_plots}
	
\end{figure}	

\begin{figure}[htbp]
	\centering
			\hspace{-2cm}
		\includegraphics[width=90mm]{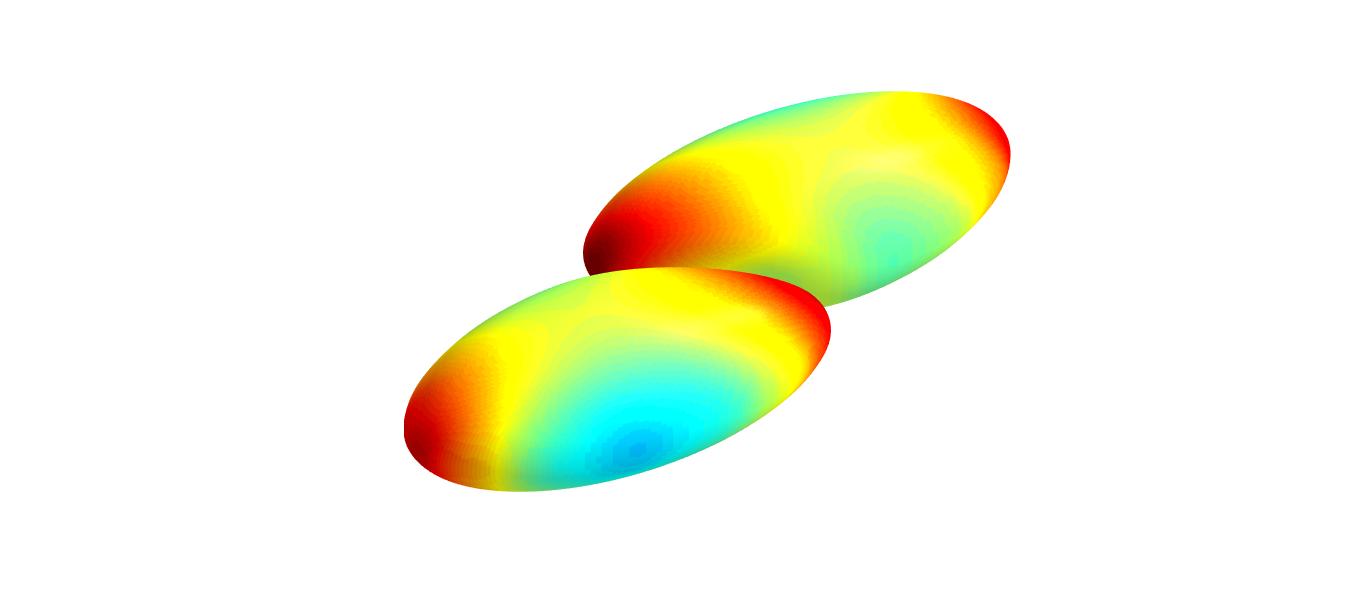}
		\quad
		\hspace{-3cm}
		\includegraphics[width=90mm]{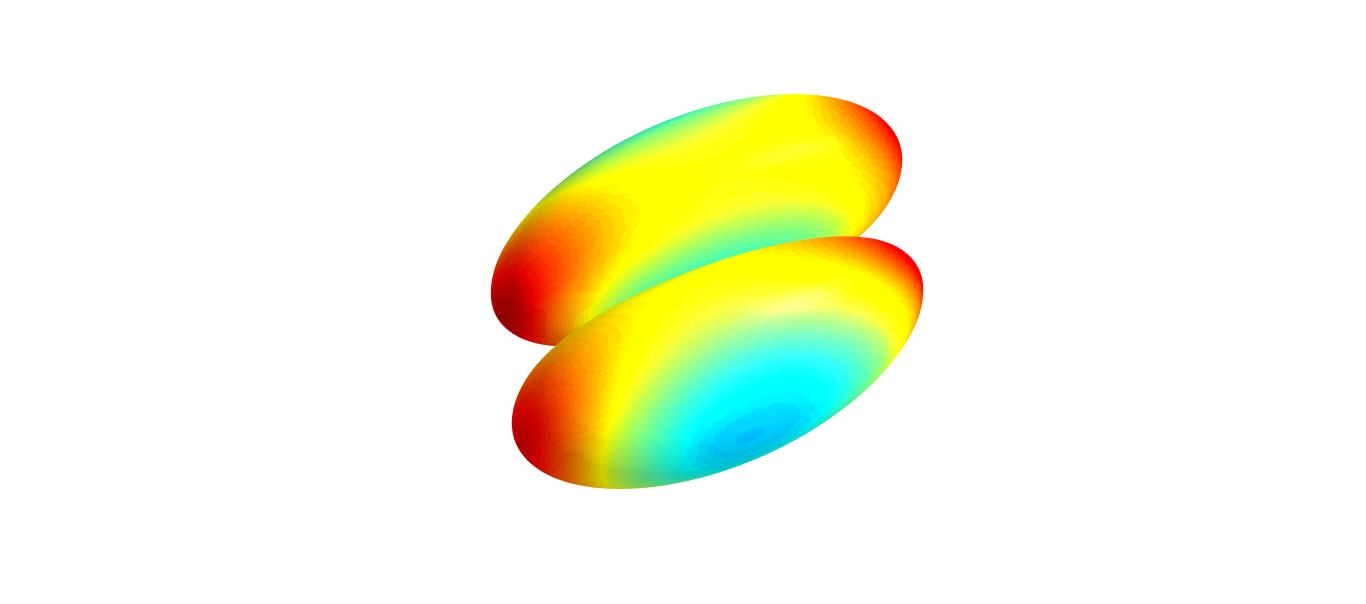}
	\caption{Two surfactant-covered drops in a shear flow at time T=18. Comparing two different viscosity ratio: $\lambda=1$ (left) and $\lambda=2$ (right).}
	\label{fig:3dropsshear1012}
	
\end{figure}

\begin{figure}[htbp]
	\centering
	\vspace{2cm}
		\includegraphics[width=64mm]{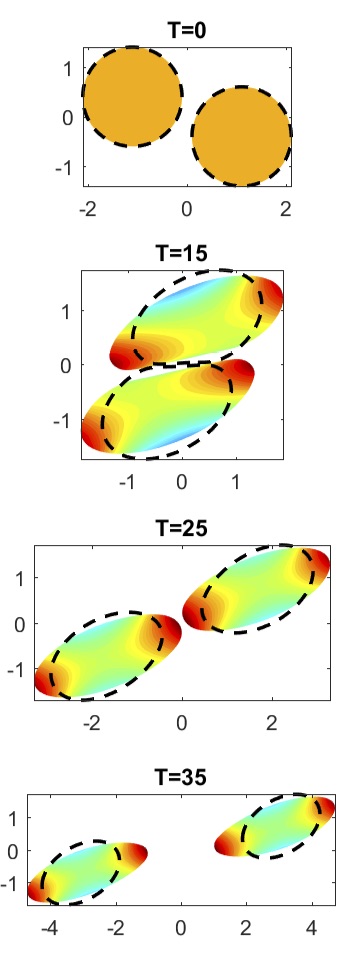}
		\quad
		\vspace{-1cm}
		\includegraphics[width=74mm]{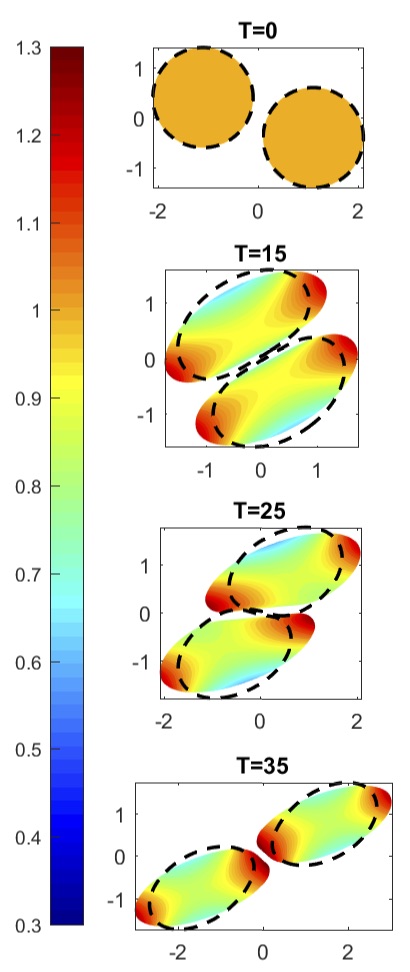}
		\vspace{1cm}
	\caption{Projections of two surfactant-covered drops interacting in a shear flow; the colors denote the surfactant concentration, the dotted line denote a clean drop. Comparing two different viscosity ratio: $\lambda=1$ (left column) and $\lambda=2$ (right column).}
	\label{fig:2dropsshear1012}
	
\end{figure}

\begin{figure}[htbp]
	\centering
	\vspace{2cm}
		\includegraphics[width=70mm]{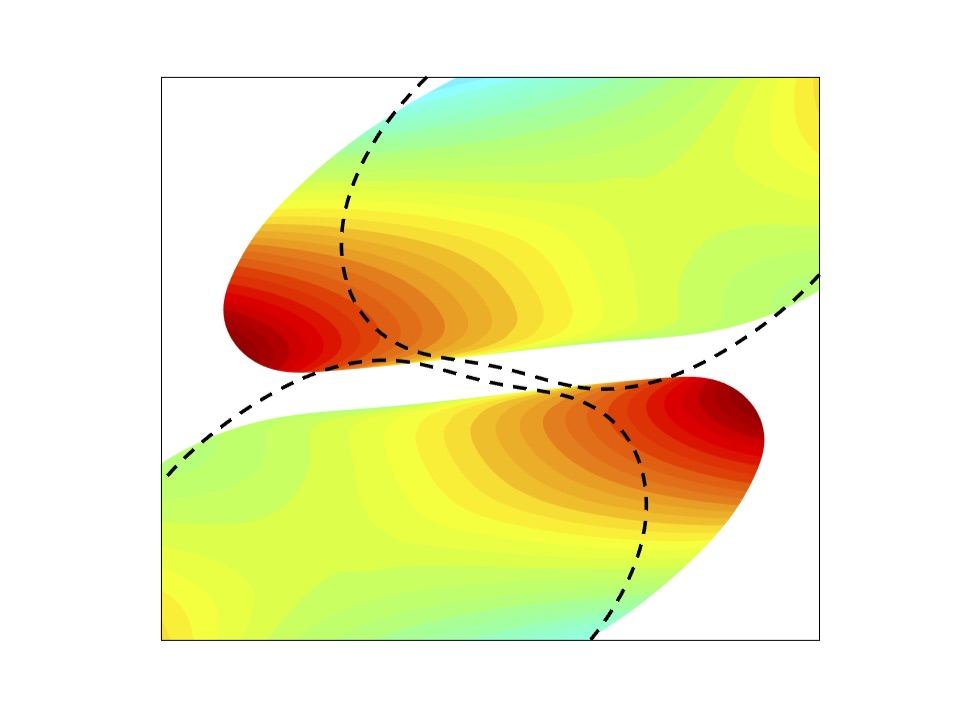}
	\caption{Zoom of the interaction between two drops for $\lambda=2$, T=25. The colors denote the surfactant concentration, the dotted line denote the contour of a clean drop }
	\label{fig:2dropsshear_zoom}
	
\end{figure}

Finally, with confidence from the successful comparisons with numerical and experimental results, we proceed with a numerical simulation of two interacting drops. Unfortunately there are no analytical, experimental or numerical results to validate the multiple surfactant covered drops experiment, but we show how our method can handle also very close interactions, how the deformation is affected by the presence of the surfactant and by the viscosity ratio $\lambda$.
In this experiment two drops are immersed in a shear flow and pushed against each other reaching a very small distance (Fig. \ref{fig:3dropsshear1012}), so that the special quadrature for the nearly-singular integrals becomes fundamental. For this experiment we set $\tilde{U}=4$ and $\tilde{L}=8$. In Fig. \ref{fig:2dropsshear1012} we show 2D projections of the drops to compare the clean drops with the surfactant-covered ones at different time-steps for two different values of $\lambda$. We can see how the surfactant acts lowering the surface tension and then making the drops deform more as compared to the clean case. The surfactant concentration is higher at the tips and this non-uniform distribution induces a strong Marangoni force. In particular we observe that for the surfactant covered drops, the minimum distance reached is never less than 0.16 for the case $\lambda=1$ and 0.15 when $\lambda=2$. When considering the clean drops, due to the smaller deformation, the minimum distance reached is instead 0.058 and 0.056 respectively for $\lambda=1,2$ (see Fig. \ref{fig:2dropsshear_zoom}) which correspond approximately to $0.25h$ that means that the special quadrature for the nearly-singular integral is absolutely needed as shown in Fig. \ref{fig:Jeffery1234}. In Fig. \ref{fig:2dropsshear1012} we compare the evolution of the drops at different times and we can observe how the higher viscosity yields a slower deformation. For these simulations the physical parameters used are the following: $\lambda=1,2$, $Ca= 0.2$, $Pe= 200$, $E=0.5$, $x_s=0.6$. We used $p=13$ combined with an adaptive upsampling and reparametrization. Note that, even using such a small number of points, we obtain an error at the final time $T=50$ which is less than 0.2\% in the volume conservation and than 0.3\% for the surfactant conservation for all cases presented. The maximum number of iterations allowed for the reparametrization is set as $i_{max}=30$, but in practice we never need so many iterations since we are allowed to use an arbitrary large pseudo-timestep as shown in Section \ref{reparametrization}.

%
 

\vspace{-0.5cm}

\section{Conclusions and future work}
\label{conclusions}
In this paper, we have presented a numerical method based on a boundary integral formulation for 3D surfactant covered drops. We have showed how the method can accurately simulate multiple drops immersed in a liquid with different viscosity ($\lambda \neq 1$) also as very close interactions occur thanks to special quadrature methods for singular and nearly singular integrals.

We presented a new reparametrization method able to maintain a good quality of the surface representation, a crucial point when dealing with strong deformations as in our case. This new procedure has no pseudo-timestep restriction and can therefore ensure fast convergence, at the same time achieving spectral accuracy of the drop volumes/areas. In addition, it provides a way to redefine the representation of any surface quantity, in our case the surfactant concentration, such that it corresponds to the new drop discretization while conserving the surface quantity (pointwise error of the surfactant) to spectral accuracy. 

Since the discretization is based on a boundary integral equation, only the surfaces of the drops must be discretized as compared to the full volume, and moreover with the representations based on a spherical harmonics expansion the number of unknowns can be kept low while still giving an excellent representation of the surfaces and of any quantities that live on the drops. This means that the system sizes for the Stokes solvers are small, but the matrices are however dense. 
The integral equation is of second kind, and yields a well conditioned system for which the condition number stays constant with an increased resolution on the drop surfaces, even though it depends on the geometry of the problem. 
This implicates that as the linear system is solved iteratively with GMRES, the number of iterations does not grow with increased resolution and the computational cost is $O(n^2)$, with $n=N(p+1)(p+2)$. 
This cost can be further reduced by using a fast summation method to evaluate the matrix-vector multiply to $O(n)$ or $O(n \log n)$. The fast multipole method \cite{fmm} is one alternative, others include the FFT based Ewald methods that has been developed for periodic \cite{DL},\cite{LAK} and recently also for free-space \cite{fast_ewald} settings. 
In the simulations presented in this paper, we have not used these methods, but it will become necessary as we want to consider more drops in one simulation. 

We studied and compared different combinations of time-stepping schemes for the drop-surfactant system, with the design criteria to minimize the number of solvers of the Stokes equations relative to the achieved accuracy, since this is the computationally most costly part. The final choice (Midpoint rule with IMEX $2^{nd}$ order) is adaptive based on the estimated accuracy for both the drop and surfactant evolution, ensuring to keep the error below a given tolerance.  This method treats the diffusive term in the surfactant equation implicitly, reducing the CFL condition to first order.  This however requires a system to be solved, which we do with GMRES. To accelerate the convergence of the iterations, we introduce a preconditioner based on analytical properties of the spherical harmonics. 

This is only one of the advantages of using a spherical harmonics representation: a Galerkin formulation (Section \ref{galerkin}) results to be cheaper than a discrete solver for the physical variables; in the reparametrization procedure the spherical harmonics representation can ensure spectral decay of the error and no pseudo-time restrictions. Also when computing singular and nearly-singular integrals, the spherical harmonics representation plays a fundamental role: in the first case (Section \ref{singular_integration}) we use analytical expressions for the rotation of the spherical harmonics coefficients (\ref{appendix2}), and when computing the interpolation for the nearly-singular quadratures (Section \ref{nearly-sing}), the on-surface value is based on a spherical harmonics interpolation. 
When two or more drops are getting close enough that the special quadrature evaluation is needed, the evaluation of the nearly singular integrals become expensive. Since no error estimates are available to guide parameter selection in the special quadrature method, choices that are unnecessarily computationally expensive might be made to ensure sufficient accuracy. 
In \cite{Ludvig_QBX} af Klinteberg and Tornberg developed a quadrature method denoted Quadrature by Expansion (QBX) \cite{Kloeckner} for the case of spheroidal rigid particles in Stokes flow. The known geometry of the rigid particles allows for precomputations that makes the method very efficient. For deformable drops, such precomputations are not possible. 
In \cite{AQBX} af Klinteberg and Tornberg developed an adaptive QBX method (AQBX) in two dimensions. Given a user specified error tolerance, parameters are selected on the fly to meet this tolerance, using error estimates derived in \cite{LaK_estimation}. In the near future, we are aiming to design an adaptive QBX method for nearly singular evaluation for deformable drops, extending the applicability of the error analysis and methodology in \cite{LaK_estimation,AQBX} to three dimensions. The goal is to minimize the cost of the integral evaluation, while maintaining errors below a set tolerance.

\section*{Acknowledgments}
\label{acknowledgments}
This work has been supported by the Knut and Alice Wallenberg Foundation under grant no KAW 2013.0339 and is gratefully acknowledged. \\
We thank Zydrunas Gimbutas and Shravan Veerapaneni for sharing their fast code for the rotation of the spherical harmonics expansions \cite{2} which has been used in the present work.

\newpage

 \appendix
 \section{Derivatives of Spherical Harmonics}
 \label{appendix1}
 For convenience, we will use the normalized associated Legendre functions which will be denoted by $\tilde{P}_n^m(x)$:
\begin{equation}
\tilde{P}_n^m(x)=\sqrt{\frac{2n+1}{(4\pi)}\frac{(n-m)!}{(n+m)!}}P_n^m(x).
\end{equation}
It follows that:
\begin{equation}
\begin{aligned}
\tilde{P}_n^{-m}(x) &=\sqrt{\frac{2n+1}{(4\pi)}\frac{(n+m)!}{(n-m)!}}P_n^{-m}(x)\\
&= \sqrt{\frac{2n+1}{(4\pi)}\frac{(n+m)!}{(n-m)!}}(-1)^m\frac{(n-m)!}{(n+m)!}P_n^m(x)\\
&=(-1)^m\sqrt{\frac{2n+1}{(4\pi)}\frac{(n-m)!}{(n+m)!}}P_n^m(x)=(-1)^m\tilde{P}_n^m(x),
\end{aligned}
\end{equation}
so that
\begin{equation}
f(\theta,\phi)= \sum_{n=0}^p \sum_{m=-n}^n f_n^m Y_n^m(\theta,\phi)=\sum_{n=0}^p \sum_{m=-n}^n f_n^m \tilde{P}_n^{m}(cos(\theta))e^{im\phi}.
\end{equation}
To compute the normalized associated Legendre functions, we can use the following recursive relations ($m \geq 0$):
\[\begin{cases}
\tilde{P}_0^0 &= \sqrt{\frac{1}{4\pi}},\\
\tilde{P}_m^m &=-c_m\sqrt{1-x^2}\tilde{P}_{m-1}^{m-1},\\
\tilde{P}_{m+1}^m &=\sqrt{2m+3}x\tilde{P}_{m}^{m},\\
\tilde{P}_{n}^m &=a_n^m(x\tilde{P}_{n-1}^{m}+b_n^m\tilde{P}_{n-2}^{m})
\end{cases}
\]
where $a_n^m =\sqrt{\frac{4n^2-1}{n^2-m^2}}$, $b_n^m=-\sqrt{\frac{(n-1)^2-m^2}{4(n-1)^2-1}}$, $c_n^m=\sqrt{1+\frac{1}{2m}}$.

The $k$-th derivatives of $f$ are then computed by:
\begin{equation}
\frac{\partial^k f}{\partial \phi^k}(\theta,\phi)=\sum_{n=0}^p \sum_{m=-n}^n (im)^k f_n^m Y_n^m(\theta,\phi)
\end{equation}
\begin{equation}
\frac{\partial^k f}{\partial \theta^k}(\theta,\phi)=\sum_{n=0}^p\sum_{m=-n}^n f_n^m \frac{\partial^k}{\partial \theta^k}\tilde{P}_n^{m}(cos \theta) e^{im\phi}
\end{equation}


where we apply recursion formulas for computing the Associated Legendre Functions derivatives. For $k=1,2$ we obtain:
\begin{equation}
\begin{aligned}
\frac{d}{d\theta}\tilde{P}_n^m(\cos \theta)&=\alpha_n^m \tilde{P}_n^{m+1}(\cos \theta)+m\cot \phi \tilde{P}_n^m(\cos \theta)\\
&= -\sqrt{d_n^m}\tilde{P}_n^{m-1}(\cos \theta)-m\cot \theta \tilde{P}_n^m(\cos \theta)
\end{aligned}
\end{equation}
\begin{equation}
\begin{aligned}
\frac{d^2}{d\theta^2}\tilde{P}_n^m(\cos \theta)&=\alpha_n^m \alpha_n^{m-1}  \tilde{P}_n^{m+2}(\cos \theta)+\alpha_n^m (2m+1)\cot \theta \tilde{P}_n^{m+1}(\cos \theta)\\
& \ \ \ -m^2\tilde{P}_n^{m}(\cos \theta)\\
&=\sqrt{d_n^m}\tilde{D}_n^{m+1}(\cos \theta)-m(1+\cot^2 \theta)\tilde{P}_n^{m}(\cos \theta)\\
& \ \ \ +m\cot \theta \tilde{D}_n^{m}(\cos \theta)\\
&= \sqrt{d_n^m}\cot \theta \tilde{P}_n^{m-1}(\cos \theta)\\
& \ \ \ +(-d_n^m+m+m(1+m)\cot^2 \theta) \tilde{P}_n^m(\cos \theta)
\end{aligned}
\end{equation}

where $d_n =(n+m)(n-m+1)$, $\alpha_n^m =\sqrt{(n-m)(n+m+1)}$.

\section{Rotation of Spherical Harmonics}
 \label{appendix2}
In order to compute singular integrals, we need to rotate the pole of the spherical harmonics expansion into a general grid location $(\beta,\alpha)$. In the rotated system, the function $f$ defined in (\ref{eq:spharm_ex}) can be expressed as (\cite{2},\cite{3},\cite{4}):
\begin{equation}
f(\theta',\phi')=\sum_{n=0}^p \sum_{m'=-n}^n f_n^{m'} (\alpha, \beta, \gamma) Y_n^{m'}(\theta',\phi')
\end{equation}
where $ (\theta',\phi') $ denote the new coordinates,
\begin{equation}
f_n^{m'} (\alpha, \beta, \gamma)= \sum_{m'=-n}^n D_n^{m',m} (\alpha, \beta, \gamma) f_n^{m}
\end{equation}
and the standard Euler angles $ (\alpha, \beta, \gamma) $ define the rotation using the $z-y-z$ convention, meaning that we first rotate by an angle $\alpha$ respect to the $z$-axis, then by an angle $\beta$ respect to the new $y$-axis and then again by an angle $\gamma$ respect to the new $z$-axis. The coefficients of the transformation are given by \cite{5}:
\begin{equation}
D_n^{m',m} (\alpha, \beta, \gamma) = e^{im\gamma}d_n^{m',m}(\beta)e^{im\alpha},
\end{equation}
where
\begin{equation}
d_n^{m',m} (\beta) = (-1)^{m'-m}[(n+m')!(n-m')!(n+m)!(n-m)!]^{1/2} S_n^{m',m}
\end{equation}
with
\begin{equation}
S_n^{m',m} (\beta) = \sum_{s=max(0,m-m')}^{min(n+m,n-m')}(-1)^s \frac{(cos \frac{\beta}{2})^{2(n-s)+m-m'} (sin \frac{\beta}{2})^{2s-m+m'}}{(n+m-s)!s!(m'-m+s)!(n-m'-s)!}.
\end{equation}


\clearpage
\newpage


\section*{References}
\begin{thebibliography}{00}
\bibitem{fast_ewald} L. af Klinteberg, D. Saffar Shamshirgar, A.-K. Tornberg. 
\emph{Fast Ewald summation for free-space Stokes potentials,} 
Res. Math. Sci., 4:1, 2017

\bibitem{AQBX}  L. af Klinteberg, A.-K. Tornberg.
\emph{Adaptive quadrature by expansion for layer potential evaluation in two dimensions,} arXiv:1704.02219 [math.NA], 2017

\bibitem{Ludvig_QBX}  L. af Klinteberg, A.-K. Tornberg.
\emph{A fast integral equation method for solid particles in viscous flow using quadrature by expansion,} J. Comput. Phys., 326:420-445, 2016

\bibitem{LaK_estimation} L. af Klinteberg, A.-K. Tornberg. 
\emph{Error estimation for quadrature by expansion in layer potential evaluation,} 
Adv. Comput. Math., 43(1):195-234, 2017

\bibitem{LAK} L. af Klinteberg, A.-K. Tornberg. 
\emph{Fast Ewald summation for Stokesian particle suspensions,} 
Int. J. Numer. Meth. Fluids, 76(10):669-698, 2014

\bibitem{IMEX}  U.M. Ascher, S.J. Ruuth, E.J. Spiteri.
\emph{Implicit-explicit Runge-Kutta methods for time-dependent partial differential equations,} 
Appl. Numer. Math., 25(2-3):151-167, 1997

\bibitem{SPHARM} K. Atkinson, W. Han.
\emph{Spherical harmonics and approximations on the unit sphere: an introduction,} 
Lecture Notes in Mathematics 2044, Springer, 2012

\bibitem{18}  I. B. Bazhlekov, P. D. Anderson, H. E.H. Meijer.
\emph{Numerical investigation of the effect of insoluble surfactants on drop deformation and breakup in simple shear flow,} 
J. Colloid Interf. Sci., 298:269-394, 2006

\bibitem{12} B.J. Bentley, L.G. Leal.
\emph{An experimental investigation of drop deformation and breakup in steady, two-dimensional linear flows,}
J. Fluid Mech., 167:241-283, 1986

\bibitem{19}  J.P. Boyd.
\emph{Chebyshev and Fourier spectral methods,} 
Dover Publications Inc., 1999

\bibitem{Bremer2012} J. Bremer, Z. Gimbutas.
\emph{A Nystr\"om method for weakly singular integral operators on surfaces,}
J. Comput. Phys., 231(14):4885-4903, 2012

\bibitem{Bremer2013} J. Bremer, Z. Gimbutas.
\emph{On the numerical evaluation of the singular integrals of scattering theory,}
J. Comput. Phys., 251:327-343, 2013

\bibitem{4} G. Burel, H. Henocq.
\emph{Determination of the orientation of 3D objects using spherical harmonics,}
Graph. Model. Im. Proc., 57(5):400-408, 1995

\bibitem{15} A.T. Chwang, T.Y.-T. Wu.
\emph{Hydromechanics of low-Reynolds-number flow. Part 2. Singularity method for Stokes flows,} 
J. Fluid Mech., 67(4):787-815, 1975

\bibitem{16}  J.E. Dennis, Jr., R.B. Schnabel.
\emph{Numerical methods for unconstrained optimization and nonlinear equations,} 
Classics in Applied Mathematics SIAM, 1996

\bibitem{11} P. Dimitrakopoulus.
\emph{Interfacial dynamics in Stokes flow via a three-dimensional fully-implicit interfacial spectral boundary element algorithm,}
J. Comput. Phys., 225:408-426, 2007 

\bibitem{feigl} K. Feigl, D. Megias-Alguacil, P. Fischer, E. J. Windhab.
\emph{Simulation and experiments of droplet deformation and orientation in simple shear flow with surfactant,}
Chem. Eng. Sci., 62:3242-3258, 2007

\bibitem{8} M. Ganesh, I.G. Graham.
\emph{A high-order algorithm for obstacle scattering in three dimensions,}
J. Comput. Phys., 198:211-242, 2004

\bibitem{3} Z. Gimbutas, L. Greengard.
\emph{A fast and stable method for rotating spherical harmonic expansions,}
J. Comput. Phys., 228:5621-5627, 2009

\bibitem{2} Z. Gimbutas, S. K. Veerapaneni.
\emph{A fast algorithm for spherical grid rotations and its application to singular quadrature,}
SIAM J. Sci. Comput., 35(6):A2738-A2751, 2013

\bibitem{9} I.G. Graham, I.H. Sloan.
\emph{Fully discrete spectral boundary integral methods for Helmholtz problems on smooth closed surfaces in $\Re^3$,}
Numer. Math., 92:289-323, 2002 

\bibitem{cut_zahedi} P. Hansbo, M.G. Larson, S. Zahedi.
\emph{A cut finite element method for a Stokes interface problem,}
Appl. Numer. Math., 85:90-114, 2014

\bibitem{22} Y.T. Hu, A. Lips.
\emph{Estimating surfactant surface coverage and decomposing its effect on drop deformation,} 
Phys. Rev. Lett., 91, 2003

\bibitem{14}  G.B. Jeffery.
\emph{The motion of ellipsoidal particles immersed in a viscous fluid,} 
Proc. R. Soc. A Math. Phys. Eng. Sci., 102(715):161-179, 1922

\bibitem{surf1} H.N. Joensson, M.L. Samuels, E.R. Brouzes, M. Medkova, M. Uhle, D.R. Link, H.A. Svahn.
\emph{Detection and analysis of cell surface biomarkers,} Angew. Chem. Int. Ed. Engl., 48(14):2518-21, 2009

\bibitem{surf_app_SRT} A.H. Jobe.
\emph{Drug therapy: pulmonary surfactant therapy,}
N. Engl. J. Med., 328:861-868, 1993

\bibitem{AK_2D_surfactant}  S. Khatri, A.-K. Tornberg.
\emph{A numerical method for two phase flows with insoluble surfactants,}
Comput. Fluids, 49(1):150-165, 2011 

\bibitem{AK-2D-soluble} S. Khatri, A.-K. Tornberg.
\emph{An embedded boundary method for soluble surfactants with interface tracking for two-phase flows,}
J. Comput. Phys., 256:768-790, 2014

\bibitem{Kloeckner} A. Kl\"ockner, A. Barnett, L. Greengard, M. O'Neil.
\emph{Quadrature by expansion: a new method for the evaluation of layer potentials,} 
J. Comput. Phys., 252:332-349, 2013

\bibitem{libro_leal} L.G. Leal.
\emph{Advanced transport phenomena - Fluid mechanics and convective transport processes,} 
Cambridge University Press, 2006

\bibitem{lipoz} X. Li, C. Pozrikidis.
\emph{The effect of surfactants on drop deformation and on the rheology of dilute emulsions in Stokes flow,}
J. Fluid Mech., 341:165-194, 1997

\bibitem{DL} D. Lindbo, A.-K. Tornberg. 
\emph{Spectrally accurate fast summation for periodic Stokes potentials,} 
J. Comput. Phys., 229(23):8994-9010, 2010

\bibitem{SDC1}  M. Minion.
\emph{Semi-implicit spectral deferred correction methods for ordinary differential equations,} 
Comm. Math. Sci., 1(3):471-500, 2003

\bibitem{fast_spharm1}  M.J. Mohlenkamp.
\emph{A fast transform for spherical harmonics,} 
J. Fourier anal. appl., 5(2):159-184, 1999

\bibitem{Muradoglu}  M. Muradoglu, G. Tryggvason.
\emph{A front-tracking method for computation of interfacial
flows with soluble surfactants,}
J. Comput. Phys., 227:2238-2262, 2008

\bibitem{orszag}  S.A. Orszag.
\emph{Fourier series on spheres,}
Mon. Weather Rev., 102:56-75, 1974

\bibitem{marangoni} Y. Pawar, K.J. Stebe.
\emph{Marangoni effects on drop deformation in an extensional flow: the role of surfactant
physical chemistry. I. Insoluble surfactants,}
Phys. Fluids, 8:1738-1751, 1996

\bibitem{10} C. Pozrikidis.
\emph{Boundary integral and singularity methods for linearized viscous flow,}
Cambridge Texts in Applied Mathematics, 1992
https://www.facebook.com/
\bibitem{SDC2}  B. Quaife, G. Biros.
\emph{Adaptive time stepping for vesicle suspensions,} 
J. Comput. Phys., 306:478-499, 2016

\bibitem{1} A. Rahimian, S. K. Veerapaneni, D. Zorin, G. Biros.
\emph{Boundary integral method for the flow of vesicles with viscosity contrast in three dimensions,}
J. Comp. Phys., 298:766-786, 2015

\bibitem{7} N. Schaeffer.
\emph{Efficient spherical harmonic transforms aimed at pseudospectral numerical simulations,}
Geochem. Geophys. Geosyst., 14:751-758, 2013

\bibitem{microfluidics_new} N. Shembekar, C. Chaipan, R. Utharala, C. A. Merten.
\emph{Droplet-based microfluidics in drug discovery, transcriptomics and high-throughput molecular genetics,}
Lab Chip, 16:1314, 2016

\bibitem{surf2} Y. Skhiri et al., 
\emph{Dynamics of molecular transport by surfactants in emulsions,}
Soft Matter, 8:10618-10627, 2012

\bibitem{17}  H.A. Stone.
\emph{A simple derivation of the time-dependent convective-diffusion equation for surfactant transport along a deforming interface,} 
Phys. Fluids A: Fluid Dynamics, 2:111, 1990

\bibitem{drops_review} H.A. Stone.
\emph{Dynamics of drop deformation and breakup in viscous fluids,}
Annu. Rev. Fluid Mech., 26:65, 1994

\bibitem{Stone2D} H.A. Stone, L.G. Leal.
\emph{The effects of surfactants on drop deformation and breakup,}
J. Fluid Mech., 220:161-186, 1990

\bibitem{microfluidics} S.Y. Teh, R. Lin, L.H. Hungand, A.P. Lee.
\emph{Droplet microfluidics,}
LabChip., 8(2):198-220, 2008

\bibitem{electric_surfactant} K.E. Teigen, K.Y. Lervag, S. T. Munkejord.
\emph{Sharp interface simulations of surfactant-covered drops in electric fields,}
V European Conference on Computational Fluid Dynamics ECCOMAS CFD, 2010

\bibitem{soluble_3D} K.E. Teigen, P. Song, J. Lowengrub, A. Voigt.
\emph{A diffuse-interface method for two-phase flows with soluble surfactants,}
J. Comput. Phys., 230:375-393, 2011

\bibitem{20}  S. Tlupova and J.T. Beale.
\emph{Nearly singular integrals in 3D Stokes flow,}
Commun. Comput. Phys., 14(5):1207-1227, 2013

\bibitem{fmm} A.-K. Tornberg, L. Greengard. 
\emph{A fast multipole method for the three-dimensional Stokes equations,} 
J. Comput. Phys., 227:1613-1619, 2008

\bibitem{6} S.K. Veerapaneni, A. Rahimian, G. Biros, D. Zorin.
\emph{A fast algorithm for simulating vesicle flows in three dimensions,}
J. Comp. Phys., 230:5610-5634, 2011

\bibitem{5} E.P. Wigner.
\emph{Group theory and its application to the quantum mechanics of atomic spectra,}
Academic Press, 1959

\bibitem{summary_microfluidics} M. W\"orner.
\emph{Numerical modeling of multiphase flows in microfluidics and micro process engineering: a review of methods and applications,}
Microfluid Nanofluid, 12:841-886, 2012

\bibitem{Xu-2D-soluble} K. Xu, M.R. Booty, M. Siegel.
\emph{Analytical and computational methods for two-phase flow with soluble surfactant,}
SIAM J. Appl. Math., 73(1):523-548, 2013

\bibitem{Xu2D} J.-J. Xu, Z. Li, J. Lowengrub, H. Zhao.
\emph{A level-set method for interfacial flows with surfactant,}
J. Comput. Phys., 212(2):590-616, 2006

\bibitem{13} L. Ying, G. Biros, D. Zorin.
\emph{A high-order 3D boundary integral equation solver for elliptic PDEs in smooth domains,}
J. Comput. Phys., 219:247-275, 2005

\bibitem{yonpoz}   S. Yon, C. Pozrikidis.
\emph{A finite-volume/boundary-element method for flow past interfaces in the presence of surfactants, with application to shear flow past a viscous drop,}
Comput. Fluids, 27(8):879-902, 1998

\bibitem{drops_bubble_book} Z. Zapryanov, S. Tabakova.
\emph{Dynamics of bubbles, drops and rigid particles,}
Springer Science+Business Media Dordrecht, 1999

\bibitem{21}  H. Zhao, A.H. Isfahani, L.N. Olson, J.B. Freund.
\emph{A spectral boundary integral method for flowing blood cells,} 
J. Comput. Phys., 229(10):3726-3744, 2010

\bibitem{zinchenko} A.Z. Zinchenko, M.A. Rother, R.H. Davis.
\emph{A novel boundary-integral algorithm for viscous interaction of deformable drops,} 
Phys. Fluids, 9:1493, 1997

\end{thebibliography}
\end{document}